\documentclass[10pt]{amsart}
\usepackage{amsfonts}
\usepackage{amsmath}
\usepackage{amsthm}
\usepackage{amssymb}
\usepackage{latexsym}
\usepackage{multicol}
\usepackage{verbatim}
\input xy
\xyoption{all}

\advance\textwidth by 1.2in \advance\oddsidemargin by -.6in
\advance\evensidemargin by -.6in
\newtheorem{cor}{Corollary}[section]
\newtheorem{lem}[cor]{Lemma}

\newtheorem{prop}[cor]{Proposition}

\newcommand\ep{\epsilon}

\newcommand{\ext}{\text{Ext}^1_{L^\sigma(\lie g)}}

\theoremstyle{definition}
\newtheorem{defn}[cor]{Definition}
\theoremstyle{definition}
\newtheorem{thm}{Theorem}

\newtheorem*{rem}{Remarks}

\newenvironment{pf}{\proof}{\endproof}
\newcounter{cnt}
\newenvironment{enumerit}{\begin{list}{{\hfill\rm(\roman{cnt})\hfill}}{%
\settowidth{\labelwidth}{{\rm(iv)}}\leftmargin=\labelwidth%
\advance\leftmargin by
\labelsep\rightmargin=0pt\usecounter{cnt}}}{\end{list}}

\theoremstyle{remark}

\numberwithin{equation}{section} 

\begin{document}

\newcommand{\thmref}[1]{Theorem~\ref{#1}}
\newcommand{\secref}[1]{Section~\ref{#1}}
\newcommand{\lemref}[1]{Lemma~\ref{#1}}
\newcommand{\propref}[1]{Proposition~\ref{#1}}
\newcommand{\corref}[1]{Corollary~\ref{#1}}
\newcommand{\remref}[1]{Remark~\ref{#1}}
\newcommand{\defref}[1]{Definition~\ref{#1}}
\newcommand{\er}[1]{(\ref{#1})}
\newcommand{\id}{\operatorname{id}}
\newcommand{\tensor}{\otimes}
\newcommand{\nc}{\newcommand}
\newcommand{\rnc}{\renewcommand}
\newcommand{\qbinom}[2]{\genfrac[]{0pt}0{#1}{#2}}
\nc{\cal}{\mathcal} \nc{\goth}{\mathfrak} \rnc{\bold}{\mathbf}
\renewcommand{\frak}{\mathfrak}
\newcommand{\desc}{\operatorname{desc}}
\newcommand{\Maj}{\operatorname{Maj}}
\renewcommand{\Bbb}{\mathbb}
\nc\bpi{{\mbox{\boldmath $\pi$}}}
\newcommand{\lie}[1]{\mathfrak{#1}}
\makeatletter
\def\section{\def\@secnumfont{\mdseries}\@startsection{section}{1}%
  \z@{.7\linespacing\@plus\linespacing}{.5\linespacing}%
  {\normalfont\scshape\centering}}
\def\subsection{\def\@secnumfont{\bfseries}\@startsection{subsection}{2}%
  {\parindent}{.5\linespacing\@plus.7\linespacing}{-.5em}%
  {\normalfont\bfseries}}
\makeatother
\def\subl#1{\subsection{}\label{#1}}

\nc{\Cal}{\cal} \nc{\Xp}[1]{X^+(#1)} \nc{\Xm}[1]{X^-(#1)}
\nc{\on}{\operatorname} \nc{\ch}{\mbox{ch}} \nc{\Z}{{\bold Z}}
\nc{\J}{{\cal J}} \nc{\C}{{\bold C}} \nc{\Q}{{\bold Q}}
\renewcommand{\P}{{\cal P}}
\nc{\N}{{\Bbb N}} \nc\boa{\bold a} \nc\bob{\bold b} \nc\boc{\bold
c} \nc\bod{\bold d} \nc\boe{\bold e} \nc\bof{\bold f}
\nc\bog{\bold g} \nc\boh{\bold h} \nc\boi{\bold i} \nc\boj{\bold
j} \nc\bok{\bold k} \nc\bol{\bold l} \nc\bom{\bold m}
\nc\bon{\bold n} \nc\boo{\bold o} \nc\bop{\bold p} \nc\boq{\bold
q} \nc\bor{\bold r} \nc\bos{\bold s} \nc\bou{\bold u}
\nc\bov{\bold v} \nc\bow{\bold w} \nc\boz{\bold z}

\nc\ba{\bold A} \nc\bb{\bold B} \nc\bc{\bold C} \nc\bd{\bold D}
\nc\be{\bold E} \nc\bg{\bold G} \nc\bh{\bold H} \nc\bi{\bold I}
\nc\bj{\bold J} \nc\bk{\bold K} \nc\bl{\bold L} \nc\bm{\bold M}
\nc\bn{\bold N} \nc\bo{\bold O} \nc\bp{\bold P} \nc\bq{\bold Q}
\nc\br{\bold R} \nc\bs{\bold S} \nc\bt{\bold T} \nc\bu{\bold U}
\nc\bv{\bold V} \nc\bw{\bold W} \nc\bz{\bold Z} \nc\bx{\bold X}

\title[Block decomposition of twisted loop representations]{The block decomposition of finite--dimensional representations of twisted loop 
algebras}

\author{Prasad Senesi}
\address{Department of Mathematics and Statistics, University of Ottawa, Ottawa, ON K1N 6N5.} \email{jsenesi@uottawa.ca}

\maketitle

\section*{abstract}
Let  $L^\sigma(\lie g)$ be the twisted loop algebra of a simple complex Lie algebra $\lie g$ with non--trivial diagram automorphism $\sigma$.  Although the category $\cal F^\sigma$ of finite--dimensional representations of $L^\sigma(\lie g)$ is not semisimple, it can be written as a sum of indecomposable subcategories (the \textit{blocks} of the category).  To describe these summands, we introduce the twisted spectral characters for $L^\sigma(\lie g)$.  These are certain equivalence classes of the spectral characters defined by Chari and Moura for an untwisted loop algebra $L(\lie g)$, which were used to provide a description of the blocks of finite--dimensional representations of $L(\lie g)$.  Here we adapt this decomposition to parametrize and describe the blocks of $\cal F^\sigma$, via the twisted spectral characters.  
 
\section*{Introduction}

In this paper we study the category $\cal F^\sigma$ of finite--dimensional representations of a twisted loop algebra $L^\sigma(\lie g)$, where $\lie g$ is a simple complex Lie algebra and $\sigma$ a diagram automorphism of $\lie g$.  While there is extensive literature on the corresponding category $\cal F$ of finite--dimensional representations of the untwisted loop algebras $L(\lie g)$ (see, for example, \cite{CFS, CL, CM, CPweyl, FoL}), until recently the treatment of $\cal F^\sigma$ has been neglected (although the simple objects of the category of graded modules for $L^\sigma(\lie g)$ were described in \cite{CPint}).

The simple objects of $\cal F^\sigma$ were described in \cite{CFS}.  However, it is not a semisimple category -- there exist objects that are indecomposable but reducible.  But we can still write any object uniquely as a direct sum of indecomposables (all objects are finite--dimensional) - and thus the category $\cal F^\sigma$ has a decomposition into indecomposable abelian subcategories.  In such a decomposition, each indecomposable object will lie in a unique indecomposable abelian subcategory, although such a subcategory may contain many nonisomorphic indecomposables.  In this case, when complete reducibility is not at hand, it is natural to search for a description of the decomposition of the category.  This is a familiar and useful strategy in the BGG category $\cal O$, for example, where the blocks are parametrized by central characters of the universal enveloping algebra of $\lie g$.  When the category of representations is semisimple (as is the case, for example, for the finite--dimensional representations of $\lie g$), the blocks are parametrized by the isomorphism classes of simple objects.  

Some features of the category $\cal F^\sigma$ can be understood in terms of the corresponding category $\cal F$ of finite--dimensional representations of $L(\lie g)$.  In particular, any simple object of $\cal F^\sigma$ can be realized by restricting the action of $L(\lie g)$ on a suitable simple object of $V$ in $\cal F$ to the subalgebra $L^\sigma(\lie g)$.  The isomorphism classes of simple objects of $\cal F$ were classified in \cite{CPweyl}, and this classification was used recently in \cite{CFS} to provide the corresponding classification of simple objects in $\cal F^\sigma$.  There the relationship between the irreducibles in $\cal F$ and in $\cal F^\sigma$ is understood using the diagram automorphism $\sigma$ which is used in the construction of $L^\sigma(\lie g)$: $\sigma$ induces a folding on the monoid of Drinfeld polynomials of $\lie g$ (this is the map $\bor$ constructed in Section 2.4 below), the result of which is a monoid of polynomials which parametrizes the irreducible modules (equivalently, the twisted Weyl modules) of $L^\sigma(\lie g)$.

The blocks of the category $\cal F$ have been described as well.  For a simple complex Lie algebra $\lie g$, we will denote by $P$ the weight lattice of $\lie g$ and by $Q$ the root lattice of $\lie g$.  In \cite{CM}, it was shown that the blocks of $\cal F$ are parametrized by the spectral characters of $L(\lie g)$ --  these are finitely supported functions $\chi: \mathbb{C}^\times \rightarrow P/Q$.  The set of all such $\chi$ forms an additive monoid, denoted by $\Xi$.   The main result of this paper is to show that the methods used in \cite{CFS} to parametrize the simple objects in $\cal F$ can be extended to parametrize the blocks of $\cal F^\sigma$.  The diagram automorphism $\sigma$ is used to construct an equivalence relation on $\Xi$, and we show that the blocks of $\cal F^\sigma$ are parametrized by the corresponding equivalence classes of spectral characters.

This paper is organized as follows.  In Sections 1 and 2 we review the main results concerning the Weyl modules for the algebras $L(\lie g)$ and $L^\sigma(\lie g)$ given in \cite{CPweyl} and \cite{CFS}.  These are certain maximal finite--dimensional highest weight (in an appropriate sense, described below) modules for the loop algebras.  They are in bijective correspondence with the irreducible modules, whose classification plays an important role in the proof of the main theorem.  In Section 3 we first review the block decomposition of the category $\cal F$ by spectral characters carried out in \cite{CM}.  Then, after defining an equivalence relation $\sim_\sigma$ on the monoid $\Xi$, we show that the equivalence classes of spectral characters parametrize the blocks of $\cal F^\sigma$.  This is done in two steps: we show that every indecomposable module must have a twisted spectral character, and that any two irreducible modules sharing the same twisted spectral character must lie in the same indecomposable abelian subcategory.  
\vskip 12pt

\noindent{\bf Acknowledgments:}
These results are developed from the author's Ph.D. thesis, written under the supervision of V. Chari.  The author would like to thank Prof. Chari for her instruction and continuing guidance.  He would also like to thank A. Moura for his help in understanding the untwisted case, as well as E. Neher and A. Savage for many useful suggestions in the preparation of this manuscript.

\setcounter{section}{0}

\section{The untwisted loop algebras and the modules $W(\bpi)$}

\subsection{Preliminaries} Throughout the paper $\mathbb{C}$ (resp. $\mathbb{C}^\times$) denotes the set of complex (resp. non--zero complex) numbers, and $\mathbb{Z}$ (resp. $\mathbb{Z}_+$)
 the set of integers (resp. non--negative) integers. Given a Lie algebra $\lie a$ we denote by $\bu(\lie a)$  the universal enveloping
 algebra of $\lie a$ and by $L(\lie a)$ the loop algebra of $\lie a$. Specifically, we have
 \begin{equation*}L(\frak a) = \frak a\otimes \mathbb{C}[t,t^{-1}],\end{equation*}
with commutator given by
\begin{equation*}
 [x\otimes t^r, y\otimes t^s]
=[x,y]\otimes t^{r+s}\end{equation*} for $x,y\in\frak a$, $r,s\in\mathbb{Z}$. We identify $\lie a$ with the subalgebra $\lie a\otimes 1$ of $L(\lie
a)$. 

 Let  $\frak{g}$ be any finite-dimensional  complex simple
Lie algebra and $\lie h$ a Cartan subalgebra of $\lie g$.  Let 
$W$ be the corresponding Weyl group, and $w_0$ the longest element of $W$.  Let  $R$ be
the set of roots of $\frak g$ with respect to $\lie h$, $I$ an index set for  a set of simple roots (and hence also for
the fundamental weights), $R^+$  the set of positive
roots, $Q^+$ (resp. $P^+$) the $\mathbb{Z}_+$ span of
the simple roots (resp. fundamental weights), $\theta$
 the highest root in $R^+$, and let
\[ P = P^+ \cup  -P^+, \; \;  Q = Q^+ \cup  -Q^+.\]
$P$ contains $Q$ as a sublattice.  Let $\pi: P \rightarrow P/Q$ be the canonical projection, and define a partial order $\geq$ on $P$ by setting $\lambda \geq \mu$ if $\lambda - \mu \in Q^+$.  We will write $\lambda > \mu$ if $\lambda \geq \mu$ and $\lambda \neq \mu$.  

Given $\alpha\in R$ let $\lie g_\alpha$ be the corresponding root space, we have
$$\lie g=\lie n^-\oplus \lie h\oplus\lie n^+,\ \ \lie n^\pm=\bigoplus_{\alpha\in R^+}\lie
g_{\pm\alpha}.$$  Fix a Chevalley basis $x^\pm_\alpha, h_\alpha$,
$\alpha\in R^+$ for $\lie g$   and set
$$x^\pm_i = x^\pm_{\alpha_i},\ \ h_i = h_{\alpha_i},\ \ i\in I.$$ In
particular for $i\in I$, $$[x^+_i,x^-_i]=h_i,\ \ [h_i,x^\pm_i]=\pm 2x_i^\pm.$$

We collect here some properties of a representation of a Lie algebra on a
finite--dimensional complex vector space.  If $V$ is a representation of a complex Lie algebra $\lie a$ and $V = V_0 \supseteq V_1 \supseteq V_2 \supseteq \ldots$ is a filtration of submodules of $V$, we will refer to a quotient module $V_i / V_{i+1}$ as an $\lie a$--\textit{constituent} (or just a \textit{constituent}, if the algebra is understood) of $V$.  If each constituent of a filtration is a simple $\lie a$--module, we say that the filtration is a composition series.  Although composition series are not unique, the Jordan--H\"older Theorem guarantees that $V$ has a unique list (up to isomorphism and re--ordering) of simple constituents.  The number of such simple constituents (counting multiplicities) is the \textit{length} of the module $V$.  
\subsection{Representations of a simple complex Lie algebra}

If $\lie g$ is a simple complex Lie algebra and $V$ is a finite--dimensional representation of $\lie g$, we can write $$V=\bigoplus_{\mu\in\lie
h^*}V_\mu,\ \ V_\mu=\{v\in V: h.v=\mu(h)v\ \ \forall \ h\in\lie
h\}.$$ Set $\text{wt}(V)=\{\mu\in\lie h^*: V_\mu\ne 0\}$. It is
well--known that
$$ V_\mu \ne 0\ \ \implies \mu\in P\ \ {\text {and}}\ \
w\mu\in\text{wt}(V) \ \ \forall \ \ w\in W,$$ and that $V$ is isomorphic
to a direct sum of irreducible representations; i.e., the category of finite--dimensional representations of $\lie g$ is semisimple. The set of
isomorphism classes of irreducible finite--dimensional $\lie
g$--modules is in bijective correspondence with $P^+$.  For any $\lambda\in P^+$ let $V(\lambda)$ be an element of the corresponding isomorphism class.
Then $V(\lambda)$ is generated by an element $v_\lambda$
satisfying the relations:
\begin{equation} \lie n^+.v_\lambda=0,\ \ h.v_\lambda=\lambda(h)v_\lambda,\ \ (x_i^-)^{\lambda(h_i)+1}.v_\lambda=0.
\end{equation}
The following facts are well--known (see \cite{Bour}, for example).
\begin{prop}\label{winv} Let $V$ be a finite--dimensional representation of
$\lie g$.
\begin{enumerate}
\item[(i)] For all  $w\in W$, $\mu\in P$, we have
$\text{dim}(V_\mu)=\text{dim}(V_{w\mu})$.
\item[(ii)] Let $V(\lambda)^*$ be the representation of $\lie g$ which is  
dual
to $V(\lambda)$. Then $$V(\lambda)^*\cong V(-w_0\lambda).$$
\end{enumerate}\end{prop}

\begin{prop}\label{shaded} Let $\lambda, \mu \in P^+$, and consider $\lie g$ as a $\lie g$--module via the adjoint representation.  
\begin{enumerate}
\item[(i)]If $\text{Hom}_{\lie g}(\lie g \otimes V(\mu), V(\lambda)) \neq 0$, then $\lambda - \mu \in Q$.
\item[(ii)] [{\cite[Proposition 1.2]{CM}}]  If  $\lambda-\mu\in Q$, then there exists a sequence of weights $\mu_l \in P^+$, $l=0,\cdots,
m$, with

\begin{enumerit}\item[(i)]  $\mu_0 = \mu$,   $\mu_m=\lambda$,
  and
\item[(ii)] ${\rm Hom}_{\lie g}(\lie g\otimes V(\mu_l), V(\mu_{l+1}))\ne 0,\
\ \forall\ 0\le l\le m$.
\end{enumerit}
\end{enumerate}
\end{prop}
\begin{pf}
We give the proof of (i).  Since $\lie g$ is semisimple, we have 
\[  \text{Hom}_{\lie g} (\lie g \otimes V(\mu), V(\lambda)) \neq 0 \Rightarrow \text{Hom}_{\lie g}(V(\lambda), \lie g \otimes V(\mu)) \neq 0.\]

Let $\phi$ be a nonzero element of $\text{Hom}_{\lie g}(V(\lambda), \lie g \otimes V(\mu))$, and $v_+$ a highest weight vector in $V(\lambda)$.  Then $\phi(v_+)$ is a weight vector in $ \lie g \otimes V(\mu)$, and we must have $\phi(v_+) \neq 0$.  Therefore $\lambda = \beta + \mu - \eta$, where $\beta \in R$ and $\eta \in Q^+$, hence $\lambda - \mu \in Q$.  
\end{pf}

\subsection{The monoid $\P$}\label{monoidP}
 Let $\cal P$  be the monoid of $I$--tuples   of
polynomials $\bpi=(\pi_1,\cdots,\pi_n)$ in an indeterminate $u$
with constant term one, with multiplication being defined
component--wise. For $i\in I$, $\lambda \in P^+$, and $a\in\mathbb{C}^\times$,  define
\[ \bpi_{\lambda, a} = \left( ( 1-au)^{\lambda(h_i)} \right) \in \P. \] 
Clearly any $\bpi \in\cal P$ can be written uniquely as a
product
\begin{equation}\label{stdfactor}\bpi =\prod_{k=1}^\ell\bpi_{\lambda_i,a_i},\end{equation}
 for some
$\lambda_1,\cdots,\lambda_\ell\in  P^+ \backslash \left\{ 0 \right\}$ and distinct elements $a_1,\ldots, a_\ell\in\mathbb{C}^\times$. We will refer to the scalars $a_i$ as the \textit{coordinates} of $\bpi$, and to the factorization (\ref{stdfactor}) as the \textit{standard decomposition} of $\bpi$.  Define a map $\cal P \to P^+$ by $\bpi\mapsto\lambda_\bpi=\sum_{i\in
I}\deg(\pi_i)\omega_i.$  

\subsection{The modules $W(\bpi)$ and $V(\bpi)$}

\begin{defn} An $L(\lie g)$--module $V$ is $\ell$--highest weight (`loop--highest weight') if there exists $v_+ \in V$ such that
\begin{align*} 
& V = \bu(L(\lie g)).v_+, \\
& L(\lie n^+).v_+ = 0, \text{   and}\\ 
& L(\lie h).v_+ = \mathbb{C} v.
\end{align*}
\end{defn}

For an $\ell$--highest weight $L(\lie g)$--module $V$ and $\lambda \in \lie h^*$, we set 
\[ V_\lambda = \left\{ v \in V : h.v = \lambda(h) v \ \forall h \in \lie h \right\},\]
and 
\[ V_\lambda^+ = \left\{ v \in V_\lambda : L(\lie n^+).v =0 \right\}.\]

The $\lie g$--modules $W(\bpi)$ we now define are the \textit{Weyl modules} for $\lie g$, first introduced and studied in \cite{CPweyl}.
\begin{defn}\label{weyldef}  \textit{The Weyl Modules $W(\bpi)$}.
 Let $\bpi \in \P$ with standard decomposition $\bpi =\prod_{k=1}^\ell\bpi_{\lambda_i,a_i}$, and $J_\bpi$ the left ideal of $\bu(L(\lie g))$ generated by the elements
\[    L(\lie n^+), \ \ \ (x^-_i)^{\lambda_\bpi(h_i)+1}, \ \ \ h \otimes t^r -  \sum_{i=1}^\ell a_i^r\lambda_i(h)  \]
for all $h\otimes t^k \in L(\lie h)$ and $i \in I$.  
Then we define the left $L(\lie g)$--module $W(\bpi)$ as 
\[ W(\bpi) = \frac{\bu(L(\lie g))}{J_\bpi}.\]
Let $w_\bpi$ be the image of $1$ under the canonical projection $\bu(L(\lie g)) \rightarrow W(\bpi)$.
\end{defn}

\begin{prop}[{\cite[Proposition 2.1, Proposition 3.1]{CPweyl}}]\label{affirred} 
\mbox{}
\begin{enumerate}
\item[(i)] $W(\bpi)$ has a unique irreducible quotient, every finite--dimensional irreducible $L(\lie g)$--module occurs as such a quotient, and for $\bpi \neq \bpi'$ the irreducible quotients of $W(\bpi)$ and $W(\bpi')$ are non--isomorphic.  Therefore the isomorphism classes of simple $L(\lie g)$--modules are in bijective correspondence with $\P$.  
\item [(ii)]Given $\bpi \in \P$ with standard decomposition
$\prod_{k=1}^\ell\bpi_{\lambda_i,a_i}$, we have an
isomorphism of $L(\lie g)$--modules
\[ W(\bpi)\cong\bigotimes_{k=1}^\ell W(\bpi_{\lambda_i,a_i}).\]
\item[(iii)] Let $V$ be any finite--dimensional $\ell$--highest weight $L(\lie g)$--module generated by an element $v$ satisfying 
\[ L(\lie n^+).v = 0, \; \; \;  L(\lie h).v = \mathbb{C}v.\]
Then $V$ is a quotient of $W(\bpi)$ for some $\bpi \in \P$. 
\end{enumerate}
\end{prop}

\noindent We will denote an element of the isomorphism class of simple $L(\lie g)$--modules corresponding to $\bpi \in \P$ by $V(\bpi)$.  

\section{The Twisted Algebras $L^\sigma(\lie g)$ and the modules $W(\bpi^\sigma)$}

\subsection{}\label{twistdef}  Here we review the construction of the twisted loop algebra $L^\sigma(\lie g)$ (for further details see \cite{Kac}).  This begins with a \textit{diagram automorphism} $\sigma$ of $\lie g$: a Lie algebra automorphism induced by a bijection $\sigma: I \rightarrow I$ which preserves all edge relations (and directions, where they occur) on the Dynkin diagram of $\lie g$.  One can verify by inspection that the only types for which such a non--trivial automorphism occurs are the types $A_n, D_n$ or $E_6$, and so we assume from here on that $\lie g$ is of one of these types.  Furthermore, in all types but $D_4$ there is a unique non--trivial automorphism of order 2, while for type $D_4$ there are exactly two non--trivial automorphisms (up to relabeling of the nodes of the Dynkin diagram): one of order two and one of order three.  We let $m$ be the order of $\sigma$, $G$ be the cyclic group with elements $\sigma^\epsilon$, $0 \leq \epsilon \leq m-1$, and for $i \in I$, we denote by $G_i$ the stabilizer of $i$ in $G$.  We also fix a primitive $m^{th}$ root of unity $\zeta$.  
  The automorphism $\sigma$ induces a permutation of $R$ given by $\displaystyle{\sigma: \sum_{i \in I} n_i \alpha_i \mapsto \sum_{i \in I} n_i \alpha_{\sigma(i)}}$, 
 and we have $$\sigma(\lie g_\alpha)=\lie g_{\sigma(\alpha)},\ \ \sigma(\lie h)= \lie h,\
\sigma(\lie n^\pm)=\lie n^\pm, $$  $$\lie
g=\bigoplus_{\epsilon=0}^{m-1}\lie g_\epsilon,\ \ \text{ where }\lie
g_\epsilon=\{x\in\lie g: \sigma(x)=\zeta^\epsilon x\}.$$ 
\noindent We also denote by $\sigma$ the automorphism of $\mathbb{C}^\times$ given by $\sigma: a \mapsto \zeta a$.  

Given any
subalgebra $\lie a$ of $\lie g$ which is preserved by $\sigma$,
set $\lie a_\epsilon=\lie g_\epsilon\cap\lie a$. It is known that
$\lie g_0$ is a simple Lie algebra, $\lie h_0$ is a Cartan
subalgebra of $\lie g_0$ and  that  $\lie g_\epsilon$ is an irreducible
representation of $\lie g_0$ for all $0\ \leq \epsilon\le m-1$. Moreover,
$$\lie n^\pm\cap \lie g_0=\lie n_0^\pm =\bigoplus_{\alpha\in R^+_{0}}(\lie g_0)_ {_{_{\pm \alpha}}}, $$
where we denote by $R_0$ the set of roots of the Lie algebra $\lie g_0$; the
sets $I_0$, $P_0^\pm$ etc. are defined similarly.  The set $I_0$ is in bijective correspondence with the set of $\sigma$--orbits of $I$. 

Suppose that  $\{y_i: i\in I\}$ is one of the sets $\{h_i: i\in
I\}$, $\{x_i^+:i\in I\}$ or $\{x_i^-:i\in I\}$ and assume that $i\ne n$ if $\lie g$ is of type $A_{2n}$. For $0 \leq \epsilon \leq m-1$, define subsets
$\{y_{i,\epsilon}: i\in I_0\}$ of $\lie
g_\epsilon$ by
\[ y_{i, \epsilon} = \frac{1}{\left| G_i \right|} \sum_{j=0}^{m-1} \zeta^{j \epsilon} y_{\sigma^j(i)}.\]
 
If $\lie g$ is of type $A_{2n}$, then we set,
\begin{gather*}
 h_{n,0}= 2(h_n + h_{n+1}),\ \ x^\pm_{n,0}=\sqrt{2}(x^\pm_{n}+x^\pm_{n+1}),\\
x^\pm_{n,1}= -\sqrt{2}( x^\pm_n - x^\pm_{n+1}),\ \ h_{n,1}=h_n-h_{n+1}, \\
y^\pm_{n,1} = \mp \frac14 \left[ x^\pm_{n,0}, x^\pm_{n,1} \right]
.\end{gather*} 
Then $\left\{ x^\pm_{i,0} , h_{i, 0} \right\}_{i \in I_0}$ is a Chevalley basis for $\lie g_0$; in particular $\left\{  h_{i, 0} \right\}_{i \in I_0}$ is a basis of $\lie h_0$. 
 
The subset $P^+_\sigma$ of $P_0^+$ is defined as follows\footnote[1]{When $\lie g$ is of type $A_{2n}$, the role of $\lambda$ in the representation theory of $L^\sigma(\lie g)$ is subject to an unusual constraint, described as follows.   If $V$ is some $\ell$--highest weight module generated by $v_+ \in V_\lambda^+$, the element $y^-_{n,1} \otimes t$ of $L^\sigma(\lie g)$ must act nilpotently on $v_+$.  The $\lie{sl}_2$--subalgebra corresponding to this generator is 
\[ \lie{sl}_2 \cong \left\langle y^{\pm}_{n,1} \otimes t^{\mp 1}, \frac{h_{n,0}}{2} \otimes 1 \right\rangle \subseteq L^\sigma(\lie g).\]
Therefore the usual $\lie{sl}_2$ theory requires $\lambda(\frac{h_{n,0}}{2}) \in \mathbb{Z}$.  This constraint motivates the definition of $P^+_\sigma$ given above. 
}:
\[ P^+_\sigma  = \begin{cases}  \lambda \in P_0^+ \text{ such that } \lambda(h_{n,0}) \in 2\mathbb{Z}, & \lie g \text{ of type } A_{2n}\\
 P_0^+, & \text{ otherwise,} \end{cases}\]
and we regard $\lambda\in
P_\sigma^+$ as an element of $P^+$ by
$$\lambda(h_i)=\begin{cases}
\lambda(h_{i,0}), \ \ i\in I_0,\  \mbox{ if $\lie g$\ is not of
type } A_{2n}\\
0,\ \ i\notin I_0,\\
(1-\delta_{i,n}/2)\lambda(h_{i,0}), \mbox{ if $\lie g$\ is of type }
A_{2n}.
\end{cases}
$$

Given $\lambda=\sum_{i\in I} m_i\omega_i\in P^+$
and $0\le \ep\le m-1$, define elements $\lambda(\epsilon)\in
P^+_\sigma$ by, \begin{gather*}\lambda(0)=\sum_{i\in I_0}
m_i\omega_i,\ \ \lambda(1)=\sum_{i\in I_0: \sigma(i)\ne
i}m_{\sigma(i)}\omega_i,\ \ \ {\text{if}}\ \ m=2\ \
 {\text{and}}\ \ \lie g \mbox{ not of type } A_{2n}\\
\lambda(0)=\sum_{i\in I_0} (1 + \delta_{i,n})m_i\omega_i,\ \
\lambda(1)=\sum_{i\in I_0: \sigma(i)\ne i}(1 +
\delta_{\sigma(i),n})m_{\sigma(i)}\omega_i,\ \  {\text{if}}\ \
m=2\ \
 {\text{and}}\ \ \lie g \mbox{  of type } A_{2n}\\
  \lambda(0) = m_1\omega_1 + m_2\omega_2, \ \
\lambda(1) = m_3 \omega_1, \ \ \lambda(2) = m_4 \omega_1, \ \ \
{\text{if}}\ \ m=3.\end{gather*}

\subsection{}
Let $\widetilde{\sigma}:L(\lie g) \to L(\lie g)$  be the automorphism
defined by linearly extending
\[ \widetilde{\sigma} (x\otimes t^k)= \zeta^k \sigma(x) \otimes t^k,
\]
for $x \in \frak g$, $k \in \mathbb{Z}$. Then $\tilde \sigma$ is
of order $m$ and we let $L^\sigma(\lie g)$ be the subalgebra of
fixed points of $\tilde\sigma$. Clearly,
\[L^\sigma(\lie g)\cong\bigoplus_{\ep=0}^{m-1} \lie g_\ep\otimes t^{m-\ep}\mathbb{C}[t^{m}, t^{-m}].\]
If $\lie a$ is any Lie subalgebra of $\lie g$, we set $L^\sigma(\lie a) = L(\lie a) \cap L^\sigma(\lie g)$.

\subsection{The  monoid $\P^\sigma$}  

Let $\P^\sigma$ be the monoid of $I_0$--tuples of polynomials $\bpi^\sigma = ( \pi_i )_{i \in I_0}$ in an indeterminate $u$ with constant term one, with multiplication being defined componentwise.  Let  $( \ , )$ be the form on $\lie h_0^*$ induced by
the Killing form of $\lie g_0$ normalized so that
$(\theta_0,\theta_0)=2$. Now define elements $\bpi^\sigma_{i, a}$, $\bpi^\sigma_{\lambda, a} \in \P^\sigma$ as follows: for $i\in I_0$ and $a\in\mathbb{C}^\times$,
$\lambda\in P_0^+$
 and  $\lie g$  not of type $A_{2n}$ 
$$\bpi^\sigma_{i,a}=((1-a^{(\alpha_i,\alpha_i)}u)^{\delta_{ij}}: j\in
I_0),\qquad \bpi^\sigma_{\lambda,a}=\prod_{i\in
I_0}\left(\bpi^\sigma_{i,a }\right)^{\lambda(h_i)},$$ while if
$\lie g$ is of type $A_{2n}$ we set for $i\in I_0$,
$a\in\mathbb{C}^\times$, $\lambda\in P^+_\sigma$,
$$\bpi^\sigma_{i,a}=((1-au)^{\delta_{ij}}:
j\in I_0),\qquad \bpi^\sigma_{\lambda,a}=\prod_{i\in
I_0}\left(\bpi^\sigma_{i,a
}\right)^{(1-\frac12\delta_{i,n})\lambda(h_i)}.$$

Define a  map $\cal P^+_\sigma\to
P_\sigma^+$ by $$\lambda_{\bpi^\sigma}=\sum_{i\in
I_0}(\deg\pi_i)\omega_i$$ if $\lie g$ is not  of type $A_{2n}$, and
$$\lambda_{\bpi^\sigma}=\sum_{i\in
I_0}(1+\delta_{i,n})(\deg\pi_i)\omega_i$$ if $\lie g$ is  of type
$A_{2n}$. 

It is clear that any $\bpi^\sigma \in \cal P_\sigma$ can be
written (non--uniquely) as a product
$$\bpi^\sigma=\prod_{k=1}^\ell\prod_{\ep=0}^{m-1}\bpi^\sigma_{\lambda_{k,\ep},
\zeta^\ep a_k},$$ where $\boa=(a_1,\cdots ,a_\ell)$ and $\boa^m$
have distinct coordinates; i.e., $a_i^m \neq a_j^m$ for $i \neq j$. We call any such expression a standard
decomposition of $\bpi^\sigma$.

\subsection{The map $\bor: \P \rightarrow \P^\sigma$}Given
$\bpi\in\cal P$ with standard factorization 
$$\bpi=\prod_{k=1}^\ell\bpi_{\lambda_k,a_k},
$$ 
define a map $\bor:\cal P \to\cal P^\sigma$  as follows:
$$\bor(\bpi)=\prod_{k=1}^\ell\prod_{\epsilon=0}^{m-1}\bpi^\sigma_{\lambda_k(\epsilon),\zeta^\epsilon
a_k}$$ 
(recall the definition of $\lambda_k(\epsilon)$ given in Section 2.1).  For any $\bpi \in \P$, we have $\lambda_{\bor(\bpi)} = \sum_{\ep = 0}^{m-1} \lambda_\bpi(\ep)$.
Note that $\bor$ is well defined (since the choice of the
$(\lambda_k,a_k)$ is unique) and set
$$\bor^{-1}(\bpi^\sigma)=\{\bpi\in\cal P:\bor(\bpi)=\bpi^\sigma\}.$$

\begin{lem}
\mbox{}
\begin{enumerate}
\item[(i)] Let $\lambda \in P_\sigma^+$ and $a \in \mathbb{C}^\times$.  Then
\[ \bor^{-1}(\bpi^\sigma_{\lambda, a}) = \begin{cases} \left\{ \bpi_{\lambda + \eta, a} \bpi_{-\sigma(\eta), -a} : \eta \in (P^+ - \lambda) \cap P^- \right\}, & m = 2; \\
 \left\{ \bpi_{\lambda + \eta_1 + \eta_2, a} \bpi_{-\sigma^2(\eta_1), \zeta a} \bpi_{-\sigma(\eta_2), \zeta^2 a} : \right.\\
         \ \ \ \ \ \ \ \ \left. \eta_1, \eta_2 \in P^-, \ \ (\eta_1 + \eta_2) \in P^+ - \lambda \right\}, & m = 3.  \end{cases} \]
\item[(ii)] Let $m = 2$, and  
\[ \bpi^\sigma = \prod_{i=1}^k \bpi^\sigma_{\lambda_{i,0}, a_i}\bpi^\sigma_{\lambda_{i,1}, -a_i}\]
be a standard factorization of $\bpi^\sigma \in \P^\sigma$, $\lambda_{i, \epsilon} \in P_0^+$.  Then 
\[ \bor^{-1}(\bpi^\sigma) = \prod_{i=1}^k \left\{ \bpi_{(\lambda_{i,0}) + \eta_i, a_i}\bpi_{(\lambda_{i,1}) - \sigma(\eta_i), -a_i} \left. \right| \eta_i \in (P^+ - (\lambda_{i,0}))\cap (P^- + \sigma(\lambda_{i,1})) \right\}  \]
\item [(iii)] Let $m = 3$, and 
\[ \bpi^\sigma = \prod_{i=1}^k \bpi^\sigma_{\lambda_{i,0}, a_i}\bpi^\sigma_{\lambda_{i,1}, \zeta a_i}\bpi^\sigma_{\lambda_{i,2}, \zeta^2 a_i}\]
be a standard factorization of $\bpi^\sigma \in \P^\sigma$, $\lambda_{i, \epsilon} \in P_0^+$.  Then 
\begin{align*} \bor^{-1}(\bpi^\sigma) = \prod_{i=1}^k \left\{  \right. &  \bpi_{(\lambda_{i,0}) + \eta_i + \nu_i, a_i}\bpi_{(\lambda_{i,1}) - \sigma^2(\eta_i), \zeta a_i} \bpi_{(\lambda_{i,2}) - \sigma(\nu_i), \zeta^2 a_i} :\\
 & \ \ \ \  \eta_i + \nu_i \in P^+ - (\lambda_{i,0}),\\
 & \ \ \ \ \ \ \ \ \sigma^2(\eta_i) \in P^- + (\lambda_{i,1}),\\
 & \ \ \ \ \ \ \ \ \ \ \ \ \sigma(\nu_i) \in P^- + (\lambda_{i,2}) \left. \right\},
 \end{align*}
\noindent where the product of sets written in (ii), (iii) is the set of all products. 
\end{enumerate}
\end{lem}
\begin{pf}
The statements will be proven only for $m = 3$.  The proof for the remaining cases when $m=2$ is simpler and uniform.  The proof begins with the following identity\footnote{The corresponding statement for the cases $m=2$ is the following:
\begin{center}For $m = 2$, $\lambda, \mu \in P_\sigma^+$, and $\eta \in (P^+ - \lambda)\cap (P^- + \sigma(\mu))$,  we have $
\bor(\bpi_{\lambda + \eta, a}  \bpi_{\mu-\sigma(\eta), -a} ) =\bpi^\sigma_{\lambda, a}\bpi^\sigma_{\mu, -a}
$. \end{center}  
}, whose verification is routine (recall that we regard $\lambda, \mu, \gamma \in P_\sigma^+$ as elements of $P^+$ via the embedding $I_0 \subseteq I$ as in Section 2.1):\\
\noindent For  $\lambda, \mu, \gamma \in P_\sigma^+$, and $\eta, \nu \in P$ such that 
\begin{equation*}\label{conds1} (\eta + \nu) \in P^+ - \lambda,\; \; \; \;  \sigma^2(\eta) \in P^- + \mu,  \; \; \; \;  \sigma(\nu) \in P^- + \gamma,\end{equation*}
\begin{equation}\label{boldr2}
 \bor(\bpi_{\lambda + \eta + \nu, a}  \bpi_{\mu-  \sigma^2(\eta), \zeta a}\bpi_{\gamma - \sigma(\nu), \zeta^2 a}) =\bpi^\sigma_{\lambda, a}\bpi^\sigma_{\mu, \zeta a}\bpi^\sigma_{\gamma, \zeta^2 a}.
\end{equation}
Now we will prove the identity
\[  \bor^{-1}(\bpi^\sigma_{\lambda, a}) = \left\{ \bpi_{\lambda + \eta_1 + \eta_2, a} \bpi_{-\sigma^2(\eta_1), \zeta a} \bpi_{-\sigma(\eta_2), \zeta^2 a} :  \eta_1, \eta_2 \in P^-, \ \ (\eta_1 + \eta_2) \in P^+ - \lambda \right\}\]
given in part (i) of the lemma.  The containment $\supseteq$ is immediate from identity (\ref{boldr2}) by taking $\mu = \gamma =0$ and $\eta = \eta_1, \nu = \eta_2$.   For the opposite containment, let 
\[\bpi = \prod_{k=1}^\ell \bpi_{\rho_k, a_k} \bpi_{\mu_k, \zeta a_k}\bpi_{\gamma_k, \zeta^2 a_k} \in \bor^{-1}(\bpi^\sigma_{\lambda,a}).\]
Then we must have $\rho_k = \mu_k = \gamma_k = 0$ for all $k$ such that $a_k^3 \neq a^3$, and so without loss of generality $\bpi =   \bpi_{\rho, a} \bpi_{\mu, \zeta a}\bpi_{\gamma, \zeta^2 a}$. 
Then 
\[ \bor(\bpi) = \bpi^\sigma_{\rho(0)+ \mu(2)+ \gamma(1), a}\bpi^\sigma_{\rho(1)+ \mu(0)+ \gamma(2), \zeta a}\bpi^\sigma_{\rho(2)+ \mu(1)+ \gamma(0), \zeta^2 a}.\]
The condition $\bpi \in \bor^{-1}(\bpi^\sigma_{\lambda,a})$ then forces $\rho = \lambda - \sigma(\mu) - \sigma^2(\gamma)$.  Therefore $\bpi$ is of the form \[\bpi_{\lambda + \eta_1 + \eta_2, a}\bpi_{ \sigma^2(\eta_1), \zeta a}\bpi_{ - \sigma(\eta_2), \zeta^2 a},\] where $\eta_1 = -\sigma(\mu)$ and $\eta_2 = -\sigma^2(\gamma)$, and the proof of part (i) of the lemma is complete.  

We continue with the proof of (iii).  From the description of $\bor^{-1}(\bpi^\sigma)$ given in \cite{CFS} (Lemma 3.5), it follows that $\bor^{-1}$ is multiplicative in the sense that 
\[ \bor^{-1}(\bpi^\sigma_1 \bpi^\sigma_2)=\bor^{-1}(\bpi^\sigma_1)\bor^{-1}( \bpi^\sigma_2), \]
where the product of the sets $\bor^{-1}(\bpi^\sigma_1)\bor^{-1}( \bpi^\sigma_2)$ is the set of products.  Therefore it suffices to prove (iii) for $k =1$, and the result will now follow from the following containment:
\begin{align}\label{contain}
\bor^{-1}\left( \bpi^\sigma_{\lambda, a} \right) \bor^{-1}\left( \bpi^\sigma_{\mu, \zeta a} \right)\bor^{-1}\left( \bpi^\sigma_{\gamma, \zeta^2 a} \right) \subseteq \left\{ \right. \bpi_{\lambda + \eta + \nu, a} & \bpi_{\mu - \sigma^2(\eta), \zeta a}\bpi_{\gamma - \sigma(\nu), \zeta^2 a} :  \\
 & \eta + \nu \in (P^+ - \lambda), \nonumber \\ 
 & \hspace{0.3 in}  \eta \in (P^- + \sigma^2(\mu)), \nonumber \\
 & \hspace{0.6 in} \nu \in (P^- + \sigma(\gamma)) \left.  \right\}.\nonumber
\end{align}
\noindent To prove the containment, let $\eta_i, \nu_i \in P^+$, $i = 0,1,2$, such that 
\begin{align*}
& \bpi_0 = \bpi_{\lambda + \eta_0 + \nu_0, a}\bpi_{-\sigma^2(\eta_0), \zeta a}\bpi_{-\sigma( \nu_0), \zeta^2 a} \in \bor^{-1}(\bpi^\sigma_{\lambda,a}),\\
& \bpi_1 = \bpi_{\mu + \eta_1 + \nu_1,\zeta a}\bpi_{-\sigma^1(\eta_1), \zeta^2 a}\bpi_{-\sigma( \nu_1), a} \in \bor^{-1}(\bpi^\sigma_{\mu, \zeta a}),\\
& \bpi_2 = \bpi_{\gamma + \eta_2 + \nu_2, \zeta^2 a}\bpi_{-\sigma^2(\eta_2), a}\bpi_{-\sigma( \nu_2), \zeta a} \in \bor^{-1}(\bpi^\sigma_{\gamma, \zeta^2 a}).
\end{align*}
Then 
\begin{align*}
 \bpi_0 \bpi_1 \bpi_2 & = \bpi_{\lambda + \eta_0 + \nu_0 - \sigma(\nu_1) - \sigma^2(\eta_2), a} \bpi_{\mu + \eta_1 + \nu_1 - \sigma^2(\eta_0) - \sigma(\nu_2), \zeta a} \bpi_{\gamma + \eta_2 + \nu_2 - \sigma^2(\eta_1) - \sigma(\nu_0), \zeta^2 a}\\
 & = \bpi_{\lambda + \eta' + \nu', a} \bpi_{\mu - \sigma^2(\eta'), \zeta a}\bpi_{\gamma - \sigma(\nu'), \zeta^2 a},
\end{align*}
where 
\begin{align*}
& \eta' = \eta_0 + \sigma^2(\nu_2) - \sigma(\eta_1)- \sigma(\nu_1), \\
& \nu' = \nu_0 + \sigma(\eta_1) - \sigma^2(\eta_2) - \sigma^2(\nu_2),
\end{align*}
and it is easily verified that $\lambda + \eta' + \nu', \mu - \sigma^2(\eta'), \gamma - \sigma(\nu') \in P^+$.

\noindent From the containment (\ref{contain}) we conclude that 
\begin{align*} \bor^{-1}(\bpi^\sigma_{\lambda,a}\bpi^\sigma_{\mu, \zeta a}\bpi^\sigma_{\gamma,\zeta^2 a}) \subseteq \left\{  \right.  & \bpi_{\lambda + \eta + \nu, a} \bpi_{\mu - \sigma^2(\eta), \zeta a}\bpi_{\gamma - \sigma(\nu), \zeta^2 a} :  \\
&  \eta + \nu \in (P^+ - \lambda), \  \eta \in (P^- + \sigma^2(\mu)), \ \nu \in (P^- + \sigma(\gamma)) \left. \right\}, 
\end{align*}
and part (iii) of the lemma is established.
\end{pf}
\begin{cor}\label{boldiboldr}
\mbox{}\\
 If $m = 2$ and $\bpi = \bpi_{\lambda,a}\bpi_{\mu, -a} \in \P$, then 
\[  \bor^{-1}\left( \bor(\bpi)\right) = \left\{ \bpi_{\lambda + \eta, a}\bpi_{\mu - \sigma(\eta), -a} \left. \right| \eta \in (P^+ - \lambda))\cap (P^- + \sigma(\mu)) \right\}. 
\]
If $m = 3$ and $\bpi = \bpi_{\lambda,a}\bpi_{\mu, \zeta a} \bpi_{\gamma, \zeta^2 a} \in \P$, then 
\[  \bor^{-1}\left( \bor(\bpi)\right) = \left\{ \bpi_{\lambda + \eta + \nu, a}\bpi_{\mu - \sigma^2(\eta), \zeta a} \bpi_{\gamma - \sigma(\nu), \zeta^2 a} \left. \right| (\eta + \nu) \in P^+ - \lambda,\; \;  \sigma^2(\eta) \in P^- + \mu,  \;  \;  \sigma(\nu) \in P^- + \gamma \right\}. 
\]
\end{cor}

\subsection{The modules $W(\bpi^\sigma)$, $V(\bpi^\sigma)$}\mbox{}

\begin{defn} An $L^\sigma(\lie g)$--module $V$ is $\ell$--highest weight if there exists $v_+ \in V$ such that
\begin{align*} 
& V = \bu(L^\sigma(\lie g)).v_+, \\
& L^\sigma(\lie n^+).v_+ = 0,\\ 
& L^\sigma(\lie h).v_+ = \mathbb{C} v.
\end{align*}
\end{defn}
For an $\ell$--highest weight $L^\sigma(\lie g)$--module $V$ and $\lambda \in \lie h_0^*$, we set 
\[ V_\lambda = \left\{ v \in V : h.v = \lambda(h) v \ \forall h \in \lie h_0 \right\},\]
and 
\[ V_\lambda^+ = \left\{ v \in V_\lambda : L^\sigma(\lie n^+).v =0 \right\}.\]  

\begin{defn}\label{weyldef} \textit{The Weyl modules $W(\bpi^\sigma)$} \mbox{}

Let $\bpi^\sigma \in \P^\sigma$ with a standard factorization $\bpi^\sigma =\prod_{k=1}^\ell\prod_{\ep=0}^{m-1}\bpi^\sigma_{\lambda_{k,\ep},
\zeta^\ep a_k}$.  For $\lie g$ not of type $A_{2n}$ let $J_{\bpi^\sigma}$ be the left ideal in $\bu(L^\sigma(\lie g))$ generated by the elements
\[  L^\sigma(\lie n^+),  \ \ \ (x^-_i)^{\lambda_{\bpi^\sigma}(h_i)+1}, \ \ \ (h_{i,\ep}\otimes
t^{mk-\ep}) - \sum_{j=1}^\ell\lambda_j(h_{i,0})a_j^{mk-\ep},\]
for all $h \otimes t^k \in L^\sigma(\lie h), \ i \in I_0$, and  for $\lie g$ of type $A_{2n}$ let $J_{\bpi^\sigma}$ be the ideal generated by the elements 
\[  L^\sigma(\lie n^+),  \ \ \ (x^-_i)^{\lambda_{\bpi^\sigma}(h_i)+1}, \ \ \ (h_{i,\ep}\otimes t^{mk-\ep}) - \sum_{j=1}^\ell(1-\frac
{1}{2}\delta_{i,n})\lambda_j(h_{i,\ep})a_j^{mk-\ep}
 \]
for all $h \otimes t^k \in L^\sigma(\lie h), \ i \in I_0$.  Then we define the left $L^\sigma(\lie g)$--module $W(\bpi^\sigma)$ as 
\[ W(\bpi^\sigma) = \frac{\bu(L^\sigma(\lie g))}{J_{\bpi^\sigma}}.\]
Let $w_{\bpi^\sigma}$ be the image of $1$ under the canonical projection $\bu(L^\sigma(\lie g)) \rightarrow W(\bpi^\sigma)$.
\end{defn}
Since $L^\sigma(\lie g)$ is a subalgebra of $L(\lie g)$, any $L(\lie g)$--module $V$ is an $L^\sigma(\lie g)$--module via restriction:
\[ \xymatrix{ L^\sigma(\lie g) \ar@{^{(}->}[r]& L(\lie g) \ar[r] & \text{End}(V)}\] 
We will denote this restriction to an $L^\sigma(\lie g)$--action by $V\left. \right|_{L^\sigma(\lie g)}$.  If two $L(\lie g)$--modules $V$, $W$ are isomorphic as $L^\sigma(\lie g)$--modules, we will write $ V \cong_{L^\sigma(\lie g)} W$.  We now state two propositions concerning the Weyl modules and their irreducible quotients.  Proposition \ref{equiv} is the twisted analog of Proposition \ref{affirred}, while Proposition \ref{nextprop} describes the relationship between the untwisted and twisted modules.

\begin{prop}[{\cite [Theorem 2]{CFS}}]\label{equiv}\mbox{}
\begin{enumerate}
\item[(i)] For any $\bpi^\sigma \in \P^\sigma$, $W(\bpi^\sigma)$ has a unique irreducible quotient, which we will denote by $V(\bpi^\sigma)$, and each irreducible $L^\sigma(\lie g)$--module occurs as such a quotient.
\item[(ii)] Let $V$ be any finite--dimensional $\ell$--highest weight $L^\sigma(\lie g)$--module generated by an element $v$ satisfying 
\[ L^\sigma(\lie n^+).v = 0, \; \; \;  L^\sigma(\lie g).v = \mathbb{C}v. \]
Then $V$ is a quotient of $W(\bpi^\sigma)$ for some $\bpi^\sigma \in \P^\sigma$. 
\item[(iii)] Let $\bpi^\sigma=\prod_{k=1}^\ell\prod_{\ep=0}^{m-1}\bpi^\sigma_{\lambda_{k,\ep},
\zeta^\ep a_k}$ be a standard decomposition of $\bpi^\sigma \in\P^\sigma$.  As $L^{\sigma}(\lie g)$--modules, we have
$$
W(\bpi^\sigma) \cong \bigotimes_{k=1}^{\ell}
W\left( \prod_{\ep=0}^{m-1}\bpi^\sigma_{\lambda_{k,\ep},
\zeta^\ep a_k}\right).$$
\end{enumerate}  
\end{prop}

Let $\P_{\text{Asym}}$ be the subset of $\P$ consisting of $\bpi \in \P$ such that, given the standard decomposition $\bpi = \prod \bpi_{\lambda_i, a_i}$ we have $a_i^m \neq a_j^m$ for $i \neq j$.  For any subset $\cal S$ of $\P$, let $\cal S_{\text{Asym}} = \cal S \cap \P_{\text{Asym}}$.  The role played by $\P_{\text{Asym}}$ is described in the following proposition.

\begin{prop}[{\cite[Propositions 4.1, 4.3, 4.5]{CFS}}]\label{nextprop}\mbox{}

Let $\bpi^\sigma \in \P^\sigma$ and $\bpi \in \bor^{-1}(\bpi^\sigma)_{\text{Asym}}$.
\begin{enumerate}
\item[(i)]  $ W(\bpi)\left.\right|_{L^\sigma(\lie g)} \cong W(\bpi^\sigma)$, and $V(\bpi)\left. \right|_{L^\sigma(\lie g)} \cong V(\bpi^\sigma)$.  
\item[(ii)] Denote the representations $W(\bpi)$ of $L(\lie g)$, $W(\bpi^\sigma)$ of $L^\sigma(\lie g)$ by 
\begin{align*} & \xymatrix{L(\lie g) \ar[r]^-{\phi_{\bpi}} & \operatorname{End}(W(\bpi))}  \hspace{0.4 in} \text{and}\\
 &\xymatrix{L^\sigma(\lie g) \ar[r]^-{\phi_{\bpi^\sigma}} & \operatorname{End}(W(\bpi^\sigma)),}
 \end{align*} 
respectively.  Then there exist ideals $I_\bpi \subseteq L(\lie g)$, $I_{\bpi^\sigma} \subseteq L^\sigma(\lie g)$ such that 
\begin{enumerate}
\item[(a)] the Lie algebra homomorphism $\phi_{\bpi^\sigma}$ factors through the quotient $L^\sigma(\lie g) / I_{\bpi^\sigma}$ to a representation 
\[\xymatrix{L^\sigma(\lie g) / I_{\bpi^\sigma} \ar[r]^-{\overline{\phi_{\bpi^\sigma}}} & \operatorname{End}(W(\bpi^\sigma))}. \]
\item[(b)] There exists a Lie algebra isomorphism $\lambda: L(\lie g) / I_\bpi \stackrel{\sim}{\rightarrow} L^\sigma(\lie g) / I_{\bpi^\sigma}$.  Therefore we have the following diagram of Lie algebra homomorphisms:
\begin{equation}\label{diagram}
\xymatrix{
& L(\lie g)\ar@{->>}_{p}[d] & & L^\sigma(\lie g)\ar^-{\phi_{\bpi^\sigma}}[r] \ar@{->>}^{p_\sigma}[d]  & \operatorname{End}(W(\bpi^\sigma))\\
& L(\lie g)/I_{\bpi} \ar^{\lambda}[rr]& & L^\sigma(\lie g) / I_{\bpi^\sigma} \ar@{.>}_{\overline{\phi_{\bpi^\sigma}}}[ur]
}
\end{equation}
where $p$ and $p_\sigma$ are the canonical projections. 
\item[(c)]  Let $W(\bpi^\sigma)_{L(\lie g)}$ denote the action of $L(\lie g)$ on $W(\bpi^\sigma)$ given by the composition 
\[ x \otimes t^r. w := \overline{\phi_{\bpi^\sigma}} \circ \lambda \circ p(x\otimes t^r).w\]
(as in Diagram (\ref{diagram})), where $w \in W(\bpi^\sigma)$ and $x\otimes t^r \in L(\lie g)$.  Then $W(\bpi^\sigma)_{L(\lie g)} \cong W(\bpi)$ and $V(\bpi^\sigma)_{L(\lie g)} \cong V(\bpi)$ as $L(\lie g)$--modules.
\end{enumerate}
\end{enumerate}
\end{prop}

\begin{rem}  \mbox{}

\begin{enumerate}
\item It is clear from the Diagram (\ref{diagram}) that the action of $L(\lie g)$ on $W(\bpi^\sigma)$ -- and hence the isomorphism $W(\bpi^\sigma)_{L(\lie g)} \cong W(\bpi)$ -- depends upon the isomorphism $\lambda: L(\lie g) / I_\bpi \stackrel{\sim}{\rightarrow} L^\sigma(\lie g) / I_{\bpi^\sigma}$.  So the expression $W(\bpi^\sigma)_{L(\lie g)}$ by itself is ambiguous -- $W(\bpi^\sigma)$ has, up to isomorphism, as many $L(\lie g)$--module structures (and hence is isomorphic to as many $L(\lie g)$--Weyl modules) as there are elements $\bpi \in \bor^{-1}(\bpi^\sigma)_{\text{Asym}}$ (the isomorphisms being determined by $\lambda$).  For this reason, when necessary we will write $W(\bpi^\sigma)_{L(\lie g)} \cong W(\bpi)$ to specify which $L(\lie g)$--module structure we have chosen for $W(\bpi^\sigma)$.  Several times we will speak of `fixing an $L(\lie g)$--action' on some $L^\sigma(\lie g)$--module; by this we mean making a choice of an isomorphism $\lambda: I_\bpi \rightarrow I_{\bpi^\sigma}$ such that $W(\bpi^\sigma)_{L(\lie g)}\cong W(\bpi)$. 

\item If we have the isomorphism $W(\bpi^\sigma)_{L(\lie g)}\cong W(\bpi)$, then we also have $ W(\bpi^\sigma)_{L(\lie g)} \left. \right|_{L^\sigma(\lie g)} \cong W(\bpi^\sigma)$.  This follows from the commutativity of the diagram
\[\xymatrix{
& L(\lie g)\ar@{->>}_{p}[d] \ar@{<-^)}[rr] & & L^\sigma(\lie g) \ar^-{\phi_{\bpi^\sigma}}[r] \ar@{->>}^{p_\sigma}[d]  & \operatorname{End}(W(\bpi^\sigma))\\
& L(\lie g)/I_{\bpi} \ar^{\lambda}[rr]& & L^\sigma(\lie g) / I_{\bpi^\sigma} \ar@{.>}_{\overline{\phi_{\bpi^\sigma}}}[ur]
}\]
\item If $W(\bpi^\sigma)_{L(\lie g)} \cong W(\bpi)$, then from Remark (2) above it follows that a subspace $U$ of $W(\bpi^\sigma)$ is an $L^\sigma(\lie g)$--submodule if and only if it is an $L(\lie g)$--submodule.
\end{enumerate}
\end{rem}
  
\begin{lem}[{\cite[Proposition 3.3]{CM}}]
Let $V(\bpi)$ be an irreducible $L(\lie g)$--constituent of $W(\bpi_{\lambda, a})$.  Then $\bpi = \bpi_{\mu, a}$, where $\mu \leq \lambda$.  
\end{lem}
\begin{prop} 
Let $\bpi^\sigma = \prod_{\epsilon = 0}^{m-1} \bpi^\sigma_{\lambda_\epsilon, \zeta^{\epsilon} a}$ and 
$\bpi_{\lambda, a} \in \mathbf{r}^{-1}(\bpi^\sigma)_{\text{Asym}}$.  Then any irreducible $L^\sigma(\lie g)$-constituent of $W(\bpi^\sigma)$ is isormophic to some $V(\bpi_{\mu, a})\left.\right|_{L^\sigma(\lie g)}$, where $\mu \leq \lambda$. 
\end{prop}
\begin{pf}
We fix an $L(\lie g)$--action $W(\bpi^\sigma)_{L(\lie g)} \cong W(\bpi_{\lambda, a})$.  Let $V$ be an irreducible $L^\sigma(\lie g)$-constituent of $W(\bpi^\sigma)$.  Then $V_{L(\lie g)}$ is isomorphic to an irreducible $L(\lie g)$-constituent of $W(\bpi_{\lambda, a})$.  Therefore $V_{L(\lie g)} \cong V(\bpi_{\mu, a}), \ \mu \leq \lambda$ (by the above lemma), and so $V \cong (V_{L(\lie g)})\left. \right|_{L^\sigma(\lie g)} \cong  V(\bpi_{\mu, a})\left. \right|_{L^\sigma(\lie g)}$.
\end{pf}

\section{Block Decomposition of the category $\cal F^\sigma$} 

\subsection{Block decomposition of a category}
Let $\lie a$ be any Lie algebra, and $\cal M$ the category of its finite--dimensional representations.  Then $\cal M$ is an abelian tensor category.  Any object in $M$ can be written uniquely as a direct sum of indecomposables, and we recall the following:

\begin{defn} 
\mbox{}
\begin{enumerate}
\item [(i)] Two indecomposable objects $V_1, V_2 \in \cal M$ are \textit{linked}, written $V_1 \sim V_2$, if there do not exist 
 subcategories $\cal M_1, \cal M_2$ such that $\cal M = \cal M_1 \oplus \cal M_2$ and $V_1 \in \cal M_1$, $V_2 \in \cal M_2$.  More generally, two objects $U, V \in \cal M$ are linked if every indecomposable summand of $U$ is linked to every indecomposable summand of $V$.  We will say that a single object $V$ in $\cal M$ is linked if there exists some other object $W$ such that $V \sim W$.  The relation $\sim$, when restricted to the collection of linked objects\footnote{The relation $\sim$ of linkage is symmetric and transitive, but it is not reflexive.  For example, if $W_1$ and $W_2$ are two objects in $\cal M$ which are \textit{not} linked, then $W = W_1 \oplus W_2$ is not linked to itself - in fact $W$ is linked to nothing at all.  
}, is an equivalence relation. 
\item [(ii)]  A \textit{block} of $\cal M$ is an equivalence class of linked objects.
\end{enumerate}
\end{defn}

\begin{prop}[{\cite[Proposition 1.1]{EM}}]\label{catdecomp}
The category $\cal M$ admits a unique decomposition into a direct sum of indecomposable abelian subcategories: $\cal M = \bigoplus_{\alpha \in \Lambda}\cal M_{\alpha}$. 
\end{prop}

In fact the indecomposable abelian subcategories of this decomposition consist of the equivalence classes of linked objects.  The goal of the rest of the paper is to provide a description of these blocks in terms of data relating to the Lie algebra $L^\sigma(\lie g)$.

\begin{defn} Let $U,V\in\cal M$ be indecomposable.
We say that $U$ is strongly linked to $V$ if there exist indecomposable $\lie a$--modules $U_1,\ldots ,U_\ell$, with $U_1=U$, $U_\ell=V$ and
either ${\rm{Hom}}_{\lie a}(U_k,U_{k+1})\neq 0$ or
${\rm{Hom}_{\lie a}}(U_{k+1},U_{k})\neq 0$ for all $1\le k < \ell $. We extend this to all of $\cal M$ by saying that two
modules $U$ and $V$ are strongly linked iff every indecomposable
summand of $U$ is strongly linked to every indecomposable
summand of $V$.
\end{defn}

\begin{lem}[{\cite[Lemma 2.2, Lemma 2.5]{CM}}]
\mbox{}
\begin{enumerate}
\item[(i)] Let $V_1, V_2$ be indecomposable objects in $\cal M$.  Then $V_1 \sim V_2$ if and only if they contain submodules $U_k \subseteq V_k$, $k = 1, 2$, with $U_1 \sim U_2$.
\item [(ii)] Two modules $U, V \in \cal M$ are linked if and only if they are strongly linked.
\end{enumerate}
\end{lem}
  
Let $\cal F$ (resp., $\cal F^\sigma$) be the category of finite--dimensional $L(\lie g)$--modules (resp., of finite--dimensional $L^\sigma(\lie g)$--modules).  From here on we fix a Lie algebra $\lie g$ of type $A, D$ or $E_6$, although any of the following results stated for untwisted loop algebras are true for the loop algebra $L(\lie g)$ of any simple Lie algebra.  We set  

\subsection{The blocks of the category $\cal F$}\mbox{}
\begin{defn} The monoid $\Xi$.
 
Let $\Xi$ be the set of all functions $\chi: \mathbb{C}^\times \rightarrow P/Q$ with finite support.
 Given $\lambda\in P^+$,
$a\in\mathbb{C}^\times$, let $\chi_{\lambda, a}\in \Xi$ be defined by
\begin{equation*}\chi_{\lambda,a}(z)=\delta_a(z)\overline{\lambda},\end{equation*}
where $\overline\lambda$ is the image of $\lambda$ in $P/Q$ and
$\delta_a(z)$ is the characteristic function of $a\in\mathbb{C}^\times$.

\end{defn}
\noindent Clearly $\Xi$ has the structure of an additive monoid under pointwise addition.  For $\bpi = \prod_{k=1}^\ell \bpi_{\lambda_k, a_k}\in \P$, we set 
\[ \chi_{\bpi} := \sum_{k=1}^\ell \chi_{\lambda_k, a_k}.\]
It is immediate from the definition that the map $\bpi \mapsto \chi_{\bpi}$  is a monoid homomorphism (from a multiplicative monoid to an additive monoid).  The elements of $\Xi$ are the \textit{spectral characters} of $L(\lie g)$.   

\begin{defn}We say that a module $V\in\cal F$ has spectral
character $\chi\in\Xi$ if, for every irreducible
constituent $V(\bpi_i)$ of $V$, we have $\chi_{\bpi_i} = \chi$ .  Let   $\cal F_\chi$ be the abelian
subcategory consisting of all modules $V\in\cal F$ with spectral
character $\chi$.
\end{defn}

The main theorem in \cite{CM} is the following, which describes the blocks of $\cal F$.

\begin{thm}[{\cite[Theorem 1]{CM}}]
The blocks of the category $\cal F$ are in bijective correspondence with the spectral characters $\chi \in \Xi$. In particular, we have 
\[ \cal F = \bigoplus_{\chi \in \Xi}\cal F_{\chi}, \] 
and each $\cal F_{\chi}$ is a block.
\end{thm}

\subsection{The blocks of the category $\cal F^\sigma$}
Here we will define the \textit{twisted spectral characters} of the twisted loop algebra $L^\sigma(\lie g)$.  These will be equivalence classes of spectral characters under a certain equivalence relation $\sim_\sigma$, defined below. First we need several technical results.

The relation $\bor(\bpi) = \bpi^\sigma$ will be illustrated with the diagram $\bpi \stackrel{\bor}{\longrightarrow}  \bpi^\sigma$.  We will also write $\bpi_1 \stackrel{\bor}{\longleftrightarrow}\bpi_2$ if $\bor(\bpi_1) = \bor(\bpi_2)$.  If $\chi_{\bpi_1} = \chi_{\bpi_2}$, we will write $\bpi_1 \sim_\chi \bpi_2$.  This relation $\sim_\chi$ is clearly an equivalence, and will be illustrated with the diagram  $\bpi_1 \stackrel{\chi}{\longleftrightarrow} \bpi_2$.

\begin{lem}\label{equivclass}  Let $\bpi, \bpi' \in \P$, $\bpi = \displaystyle\prod_{k=1}^\ell \bpi_{\lambda_k, a_k}$.  Then $\bpi' \sim_\chi \bpi$ if and only if $\bpi'$  is of the form 
\[ \bpi' = \prod_{k=1}^\ell \bpi_{\lambda_k + \nu_k, a_k} \tilde \bpi, \]
where $\nu_k \in Q$ such that $\lambda_k + \nu_k \in P^+$ and $\tilde \bpi \in \P$ such that $\chi_{\tilde \bpi} = 0$. 
\end{lem}
\begin{pf}
Because $\chi$ is a monoid homomorphism, it suffices to prove the lemma in the case $\ell =1$; i.e. $\bpi = \bpi_{\lambda, a}$.  If $\lambda \in Q$, then we take $\nu_1 = 0$ and $\tilde \bpi = 1$.  Now suppose $\lambda \not \in Q$ and $\chi_\bpi = \chi_{\bpi'}$.  Let us write $\bpi' = \prod_{k=1}^r \bpi_{\mu_k, b_k}\prod_{k=1}^s \bpi_{\gamma_k, c_k}$ where all $b_j$ are pairwise distinct, $\mu_k \not \in Q$ and $\gamma_k \in Q$.  Then 
\[   \chi_\bpi = \chi_{\bpi'} \ \Rightarrow  \ \delta_{a}(z) \overline{\lambda} = \sum_{k=1}^r \delta_{b_k}(z) \overline{\mu_k} .\]
Assume $r > 1$.  Evaluating this expression at $z=a$ forces $b_j = a$ for some $1 \leq j \leq r$ and $\overline{\lambda} = \overline{\mu_j}$, hence $ \mu_j - \lambda \in Q$.  Therefore $\mu_j = \lambda + \nu$ for some $\nu \in Q$.  Next, evaluating the equality at any $z = b_k$, $k \neq j$, gives us $\overline{\mu_k} = 0$, hence $\mu_k \in Q$, a contradiction.  So we must have $r=1$.  Setting $\tilde \bpi = \prod_{k=1}^s\bpi_{\gamma_k, c_k}$, we have $\chi_{\tilde \bpi} = 0$ and $\bpi_{\lambda, a} = \bpi_{\mu_j, b_j}\tilde \bpi  = \bpi_{\lambda+ \nu, a} \tilde\bpi$, as desired.
\end{pf}

\begin{defn}
We define a relation $\sim_\sigma$ on $\Xi$ as follows: $\chi_1 \sim_\sigma \chi_2$ if there exist $\bpi_i \in \P$, $i = 1, 2$, such that $\chi_i = \chi_{\bpi_i}$ and $\bor(\bpi_1) = \bor(\bpi_2)$.
\end{defn}
It is routine to show that $\chi_1 \sim_\sigma \chi_2$ if and only if 
\begin{equation}\label{sigequiv} \sum_{\epsilon = 0}^{m-1} \sigma^\epsilon \circ \chi_1 \circ \sigma^{-\epsilon} = \sum_{\epsilon = 0}^{m-1} \sigma^\epsilon \circ \chi_2 \circ \sigma^{-\epsilon},\end{equation}
(where we regard $\sigma$ as an automorphism of $\mathbb{C}^\times$ via $a \mapsto \zeta a$ and as an automorphism of $P / Q$ via $\overline{\omega_i} \mapsto \overline{\sigma(\omega_i)}$) and therefore that $\sim_\sigma$ is an equivalence relation on $\Xi$.


\begin{defn}  The twisted spectral characters $\Xi^\sigma$ of $L^\sigma(\lie g)$, denoted $\Xi^\sigma$, are the equivalence classes of $\Xi$ with respect to the equivalence relation $\sim_\sigma$: 
\[ \Xi^\sigma = \Xi / \sim_\sigma\] 
If $\bpi_\sigma \in \P^\sigma$, we define $\chi_{\bpi^\sigma} := \overline{\chi_{\bpi}}$, where $\bpi \in \bor^{-1}(\bpi^\sigma)$.  Using the relation (\ref{sigequiv}), we can see that the binary operation $\overline{\chi_1} + \overline{\chi_2} = \overline{\chi_1 + \chi_2}$ is well--defined, hence $\Xi^\sigma$ is an abelian monoid. 
\end{defn}

\begin{defn}
We say that an $L^\sigma(\lie g)$--module $V$ has spectral character $\overline{\chi}$ if, for every irreducible $L^\sigma(\lie g)$-constituent $V(\bpi^\sigma)$ of $V$, we have $\chi_{\bpi^\sigma} = \overline{\chi}$.  Let $\cal F^\sigma_{\overline{\chi}}$ be the abelian subcategory of all $L^{\sigma}(\lie g)$-modules with spectral character $\overline{\chi}$.
\end{defn}

The main result of this paper is the following theorem.

\begin{thm} The blocks of $\cal F^\sigma$ are in bijective correspondence with $\Xi^\sigma$.  In particular, 
\[ \cal F^\sigma = \bigoplus_{\overline{\chi} \in \Xi^\sigma} \cal F^\sigma_{\overline{\chi}},\]
and each $\cal F^\sigma$ is a block.
\end{thm}

The theorem follows from the next two propositions:

\begin{prop}\label{firstprop} Any two irreducible modules in $\cal F^\sigma_{\overline{\chi}}$ are linked.\end{prop}
\begin{prop}\label{secondprop}  Every indecomposable $L^\sigma(\lie g)$--module has a twisted spectral character. \end{prop}
\noindent The remainder of the paper is devoted to the proof of these propositions.
\subsection{Proof of Proposition \ref{firstprop}}

\begin{lem}\label{sqarefreerep}
Let $\bpi_{\lambda, a} \in \P$, and suppose $\bpi \in \P$ such that $\chi_{\lambda, a} = \chi_\bpi$.  Then there exists some $\tilde\bpi \in \P_{\text{Asym}}$ such that $\chi_{\tilde\bpi} = \chi_{\lambda, a}$ and $\bor(\bpi) = \bor(\tilde\bpi)$. 
\end{lem}
\begin{pf}
Since $\bpi_{\lambda, a}\sim_\chi \bpi$, $\bpi$ must be of the form 
\[ \bpi = \bpi_{\lambda + \lambda', a} \prod_{\ep = 1}^{m-1}\bpi_{\eta_\ep, \zeta^\ep a} \prod_{i=1}^\ell \prod_{\ep = 0}^{m-1} \bpi_{\nu_{\ep,i}, \zeta^\ep b_i}, \ \ \ \ \ (\text{Lemma \ref{equivclass})}\]
for $\lambda', \eta_\ep, \nu_{\ep,i} \in Q, \lambda + \lambda' \in P^+, \eta_\ep, \nu_{\ep,i} \in P^+$, and $b_i^m \neq b_j^m \neq a^m$ for all $1 \leq i \neq j \leq \ell$.  Define 
\[ \tilde\bpi = \bpi_{\lambda + \lambda' + \sum_{\ep=1}^{m-1}\sigma^{m-\ep}(\eta_\ep), a} \prod_{i=1}^\ell \bpi_{\nu_{0,i} + \sum_{\ep=1}^{m-1}\sigma^{m-\ep}(\nu_{\ep, i}), b_i}. \]
Then $\tilde\bpi \in \P_{\text{Asym}}$, $\chi_{\tilde \bpi} = \chi_{\bpi_{\lambda, a}}$ and $\bor(\bpi) = \bor(\tilde\bpi)$. 
\end{pf}

\begin{lem}\label{lem17}
\mbox{}
\begin{enumerate}
\item [(i)] Let $\bpi \in \bor^{-1}(\bpi^\sigma)_{\text{Asym}}$.  Then $V(\bpi)_{L^\sigma(\lie g)} \in \cal F^\sigma_{\overline{\chi_{\bpi}}}$.
\item [(ii)] Let $V(\bpi^\sigma) \cong V(\bpi)_{L^\sigma(\lie g)}$, $\bpi \in \bor^{-1}(\bpi^\sigma)_{\text{Asym}}$.  Then $V (\bpi^\sigma) \in \cal F^\sigma_{\overline{\chi_{\bpi}}}$.
\end{enumerate}
\end{lem}
\begin{pf}
The lemma follows directly from the definitions.  For the first, note that $\chi_{\bpi^\sigma} = \overline{\chi_{\bpi}}$, and  $V(\bpi)_{L^\sigma(\lie g)} = V (\bor(\bpi)) = V (\bpi^\sigma)$.   The second is immediate from the first.  
\end{pf}

\begin{prop}
Let $ V(\bpi^\sigma_k) \in \cal F^\sigma_{\overline{\chi_k}}$ for some $\chi_k \in \Xi$, $k = 1, 2$.  Then $ V(\bpi^\sigma_1) \otimes  V(\bpi^\sigma_2) \in  \cal F^\sigma_{\overline{\chi_1}+ \overline{\chi_2}}$.
\end{prop}

\begin{pf}
Let $ V(\bpi^\sigma_k) \in \cal F_{\overline{\chi_k}}$ for some $\chi_k \in \Xi$, $k = 1, 2$.  Choose $\bpi_i \in \bor^{-1}(\bpi^\sigma_i)_{\text{Asym}}$; therefore  $\overline{\chi_i} = \overline{\chi_{\bpi_i}}$ and $(V(\bpi_1) \otimes V(\bpi_2))_{L^\sigma(\lie g)} \cong  V(\bpi^\sigma_1) \otimes  V(\bpi^\sigma_2)$.  Fix $L(\lie g)$--actions
\begin{align*} 
& V(\bpi^\sigma_i)_{L(\lie g)} \cong V(\bpi_i)   \ \  \text{and}\\
& ( V(\bpi^\sigma_1) \otimes  V(\bpi^\sigma_2))_{L(\lie g)} \cong V(\bpi_1) \otimes V(\bpi_2).
\end{align*}
Now let $V$ be an irreducible $L^\sigma(\lie g)$--constituent of $ V(\bpi^\sigma_1) \otimes  V(\bpi^\sigma_2)$.  Then $V_{L(\lie g)}$ is some irreducible $L(\lie g)$--constituent $V(\bpi)$ of $( V(\bpi^\sigma_1) \otimes  V(\bpi^\sigma_2))_{L(\lie g)} \cong V(\bpi_1) \otimes V(\bpi_2)$.  We know from the untwisted affine case that $V(\bpi)$ has spectral character $\chi_{\bpi_1} + \chi_{\bpi_2}$ (and hence $\chi_{\bpi} = \chi_{\bpi_1} + \chi_{\bpi_2}$), and by Lemma \ref{lem17} (ii), $V$ has character $\overline{\chi_{\bpi}} = \overline{\chi_{\bpi_1} + \chi_{\bpi_2}} = \overline{\chi_{\bpi_1}} + \overline{\chi_{\bpi_2}} = \overline{\chi_1} + \overline{\chi_2} $. 
\end{pf}

\begin{cor}
For all $\chi_k \in \Xi$, $k = 1, 2$, we have 
\[ \cal F^\sigma_{\overline{\chi_1}}\otimes \cal F^\sigma_{\overline{\chi_2}} \subseteq \cal F^\sigma_{\overline{\chi_1} + \overline{\chi_2}}.\]
\end{cor}
\begin{prop}\label{spectralweyl} $W(\bpi^\sigma) \in \cal F^\sigma_{\chi_{\bpi^\sigma}}$. 
\end{prop}
\begin{pf}
In view of the above corollary, it suffices to prove the lemma when $\bpi^\sigma = \prod_{\ep =0}^{m-1} \bpi^\sigma_{\lambda_\ep, \zeta^{\ep}a}$ for $a \in \mathbb{C}^\times$, $\lambda_\ep \in P^+_\sigma$.  Let $\bpi_{\lambda, a} \in \bor^{-1}(\bpi^\sigma)_{\text{Asym}}$, so that $\chi_{\bpi^\sigma} = \overline{\chi_{\lambda, a}}$, and fix an isomorphism 
\[W(\bpi^\sigma)_{L(\lie g)} \cong W(\bpi_{\lambda, a}). \]
Now let $V = V(\bpi^\sigma_1)$ be an irreducible $L^\sigma(\lie g)$--constituent of $W(\bpi^\sigma)$.   We will show that $\chi_{\bpi^\sigma} = \chi_{\bpi^\sigma_1}$.  

$V_{L(\lie g)}$ is an irreducible $L(\lie g)$--constituent of $W(\bpi^\sigma)_{L(\lie g)} \cong W(\bpi_{\lambda, a})$.  Since $W(\bpi_{\lambda, a}) \in \cal F_{\chi_{\lambda,a}}$ (\cite{CM}, Lemma 5.1), $V_{L(\lie g)} \cong V(\bpi_1)$ for some $\bpi_1 \in \P$ such that $\chi_{\bpi_1} = \chi_{\lambda, a}$.  Since $V(\bpi_1)$ is an irreducible $L(\lie g)$--constituent of $W(\bpi_{\lambda, a})$, $V(\bpi_1)$ must be of the form $V(\bpi_{\mu, a})$ for some $\mu \leq  \lambda$ (\cite{CPweyl}, Proposition 3.3).  Therefore $V(\bpi_1^\sigma) = (V(\bpi^\sigma_1)_{L(\lie g)})_{L^\sigma(\lie g)} = V(\bpi_1)_{L^\sigma(\lie g)} = V(\bor(\bpi_1))$, hence $\bpi_1 \in \bor^{-1}(\bpi^\sigma_1)$. 

 Therefore 
\[ \chi_{\bpi^\sigma_1} = \overline{\chi_{\bpi_1}} = \overline{\chi_{\bpi}} = \chi_{\bpi^\sigma}.\]\end{pf}

The following proposition provides a strong linking between certain irreducible $L^\sigma(\lie g)$--modules.
\begin{prop}\label{machine} 
\mbox{}\\
Let $a \in \mathbb{C}^\times$, $\lambda_\ep, \mu_\ep \in P^+_\sigma$, and $\displaystyle{\lambda = \sum_{\epsilon = 0}^{m-1} \sigma^{m - \epsilon} (\lambda_{\epsilon}),  \ \mu = \sum_{\epsilon = 0}^{m-1} \sigma^{m - \epsilon} (\mu_{\epsilon})}$, so that 
\[ \bpi_{\lambda, a} \in \bor^{-1}\left(\prod_{\epsilon = 0}^{m-1} \bpi^\sigma_{\lambda_\epsilon, \zeta^{\epsilon} a}\right)_{\text{Asym}}, \ \ \ \ \bpi_{\mu, a} \in \bor^{-1}\left(\prod_{\epsilon = 0}^{m-1} \bpi^\sigma_{\mu_\epsilon, \zeta^{\epsilon} a}\right)_{\text{Asym}}.\]
Assume there exists a nonzero homomorphism $p: \lie g \otimes V(\lambda) \rightarrow V(\mu)$ of $\lie g$-modules.  The following formula defines an action of a $L^\sigma(\lie g)$-module on $V(\lambda) \oplus V(\mu)$:

\begin{equation}\label{twistaction} x \otimes t^k (v,w) = (a^k xv, a^kxw + k a^{k-1}p(x \otimes v)),
\end{equation}
where $x \in \lie g_{\overline{k}}, v \in V(\lambda)$ and $w \in V(\mu)$.  Denote this $L^\sigma(\lie g)$-module by $V(\lambda, \mu, a)$.  Then 

\[ 0 \rightarrow V\left(\prod_{\epsilon = 0}^{m-1} \bpi^\sigma_{\lambda_\epsilon, \zeta^{\epsilon} a}\right) \rightarrow V(\lambda, \mu, a) \rightarrow V\left(\prod_{\epsilon = 0}^{m-1} \bpi^\sigma_{\mu_\epsilon, \zeta^{\epsilon} a}\right) \rightarrow 0 \]
is a non-split short exact sequence of $L^\sigma(\lie g)$-modules.  If $\lambda > \mu$, then there exists a canonical surjective homomorphism of $L^\sigma(\lie g)$-modules $W(\prod_{\epsilon = 0}^{m-1} \bpi^\sigma_{\lambda_\epsilon, \zeta^{\epsilon} a}) \rightarrow V(\lambda, \mu, a)$.
\end{prop}
\begin{pf}
For brevity we will write $V\left(\prod_{\epsilon = 0}^{m-1} \bpi^\sigma_{\lambda_\epsilon, \zeta^{\epsilon} a}\right) = V(\bor(\bpi_{\lambda, a}))$, $V\left(\prod_{\epsilon = 0}^{m-1} \bpi^\sigma_{\mu_\epsilon, \zeta^{\epsilon} a}\right) = V(\bor(\bpi_{\mu, a}))$.

The verifications that formula (\ref{twistaction}) gives an $L^\sigma(\lie g)$-action, and that the sequence is exact, are routine.  To prove that the sequence is non-split, assume that 
\[ V(\lambda, \mu, a) = W_1 \oplus W_2 \]
is a non-trivial decomposition of $V(\lambda, \mu, a)$ into $L^\sigma(\lie g)$-submodules.  It is immediate from its construction that the length of $V(\lambda, \mu, a)$ is 2, with constituents $V(\bor(\bpi_{\lambda, a}))$ and $V(\bor(\bpi_{\mu, a}))$.  Therefore we can assume without loss of generality that $W_1 \cong  V(\bor(\bpi_{\lambda, a}))$.  But it is clear from the description of the action of  $L^\sigma(\lie g)$ on $V(\lambda, \mu, a)$ that $V(\bor(\bpi_{\lambda, a}))$ is not a submodule of $V(\lambda, \mu, a)$.  Therefore $V(\lambda, \mu, a)$ must be indecomposable.

Let $v_+$ be a highest-weight vector of $V(\bor(\bpi_{\lambda, a}))$.  Then $\bu(L^\sigma(\lie g)).(v_+, 0)$ must be isomorphic to  $V(\bor(\bpi_{\mu, a}))$ or $V(\lambda, \mu, a)$.   If we assume that $\lambda > \mu$, then we cannot have $\bu(L^\sigma(\lie g)).(v_+, 0) \cong V(\bor(\bpi_{\mu, a}))$ by weight considerations.  Therefore if $\lambda > \mu$, then $V(\lambda, \mu, a)$ is cyclically generated by $(v_+, 0)$.  Since this element is also highest-weight with $L^\sigma(\lie g)$-weights given by $\bor(\bpi_{\lambda, a})$, it follows that $V(\lambda, \mu, a)$ is a quotient of $W(\bor(\bpi_{\lambda, a}))$. 
\end{pf}

\begin{cor}\label{linkagemachine}
\mbox{}\\
 Let $\bpi_{\lambda, a} \in \bor^{-1}(\prod_{\epsilon = 0}^{m-1} \bpi^\sigma_{\lambda_\epsilon, \zeta^{\epsilon} a})_{\text{Asym}}$, $\bpi_{\mu, a} \in \bor^{-1}(\prod_{\epsilon = 0}^{m-1} \bpi^\sigma_{\mu_\epsilon, \zeta^{\epsilon} a})_{\text{Asym}}$ as in the proposition, and assume there exists a non-zero homomorphism $p: \lie g \otimes V(\lambda) \rightarrow V(\mu)$ of $\lie g$-modules.  Then the $L^\sigma(\lie g)$-modules $V(\prod_{\epsilon = 0}^{m-1} \bpi^\sigma_{\lambda_\epsilon, \zeta^{\epsilon} a})$ and $V(\prod_{\epsilon = 0}^{m-1} \bpi^\sigma_{\mu_\epsilon, \zeta^{\epsilon} a})$ are strongly linked.
\end{cor}

The following proposition is given in \cite{CM}.  We give a proof of the proposition here to clarify the proof of the analogous statement for the twisted case.

\begin{prop}[{\cite[Proposition 2.3]{CM}}]
Let $V(\bpi_1)$, $V(\bpi_2)$ be irreducible $L(\lie g)$--modules with $V(\bpi_i) \in \cal C_{\chi}$, $i = 1,2$.  Then $V(\bpi_1)$, $V(\bpi_2)$ are strongly linked. 
\end{prop}

\begin{pf} 
Since $V(\bpi_i) \in \cal C_{\chi}$, there exist $\lambda_i, \mu_i \in P^+$, $a_i \in \mathbb{C}$, $1 \leq i \leq \ell$ with $\lambda_i - \mu_i \in Q$, $a_i \neq a_j$ and 
\[ \bpi_1 = \prod_{i=1}^\ell \bpi_{\lambda_i, a_i}, \ \ \ \ \bpi_2 = \prod_{i=1}^\ell \bpi_{\mu_i, a_i}.\]
Let us assume for simplicity here that $\ell = 2$; the more general case is a straightforward extension of this one.  By Proposition \ref{shaded}(ii) there exist sequences of weights $\left\{ \nu_i\right\}_{i=0}^q$, $\left\{ \eta_i\right\}_{i=0}^r$  with $\nu_0 = \lambda_1$, $\nu_q = \mu_1$, $\eta_0 = \lambda_2$, $\eta_q = \mu_2$, such that 
\[ \text{Hom}_{\lie g}(\lie g \otimes V(\nu_i), V(\nu_{i+1})) \neq 0, \ \ \text{Hom}_{\lie g}(\lie g \otimes V(\eta_j), V(\eta_{j+1})) \neq 0  \]
for $0 \leq i \leq q-1, \; 0 \leq j \leq r -1$.  Fix some $i$, $1 \leq i \leq q$.  By Proposition \ref{shaded}(i) either $\nu_i \geq \nu_{i+1}$ or $\nu_{i+1} \geq \nu_{i}$.  If $\nu_i \geq \nu_{i+1}$, we can conclude (\cite{CM}, Proposition 3.4) that $V(\bpi_{\nu_i, a_1})$ and $V(\bpi_{\nu_{i+1}, a_1})$ are both irreducible constituents of some quotient $M_i$ of $W(\bpi_{\nu_i, a_1})$.

If $\nu_{i+1} \geq \nu_{i}$, we use the isomorphism 
\begin{align*}
&\text{Hom}_{\lie g}(\lie g \otimes V(\nu_i), V(\nu_{i+1}))\\
& \quad \cong \text{Hom}_{\lie g}(\lie g \otimes V(\nu_{i+1})^*, V(\nu_{i})^*)\\
&  \qquad  \cong \text{Hom}_{\lie g}(\lie g \otimes V(-w_0(\nu_{i+1})), V(-w_0(\nu_{i})))
\end{align*}
to conclude that  $V(\bpi_{\nu_i, a_1})$ and $V(\bpi_{\nu_{i+1}, a_1})$ are both irreducible constituents of some quotient $M_i$ of $W(\bpi_{\nu_{i+1}, a_1})$.  

We may now assume without loss of generality that $\nu_i \geq \nu_{i+1}$ for $0 \leq i \leq q-1$.  So for all $0 \leq i \leq q-1$, $V(\bpi_{\nu_{i}, a_1}) \otimes V(\bpi_{\lambda_2, a_2})$ and $V(\bpi_{\nu_{i+1}, a_1}) \otimes V(\bpi_{\lambda_2, a_2})$ are simple constituents of $M_i \otimes V(\bpi_{\lambda_2, a_2})$.  This module, in turn, is a quotient of $W(\bpi_{\nu_i, a_1}) \otimes W(\bpi_{\lambda_2, a_2}) \cong W(\bpi_{\nu_i, a_1}\bpi_{\lambda_2, a_2})$ hence indecomposable.  Therefore $V(\bpi_{\nu_{i}, a_1}) \otimes V(\bpi_{\lambda_2, a_2})$ and $V(\bpi_{\nu_{i+1}, a_1}) \otimes V(\bpi_{\lambda_2, a_2})$ are strongly linked, and so $V(\bpi_{\lambda_1, a_1}) \otimes V(\bpi_{\lambda_2, a_2})$ and $V(\bpi_{\mu_1, a_1}) \otimes V(\bpi_{\lambda_2, a_2})$ are strongly linked.   We show similarly that $V(\bpi_{\mu_1, a_1}) \otimes V(\bpi_{\lambda_2, a_2})$ and $V(\bpi_{\mu_1, a_1}) \otimes V(\bpi_{\mu_2, a_2})$ are strongly linked to complete the proof.
\end{pf}

\begin{prop}\label{mainprop}
Let $\left\{ \lambda_i \right\}_{i=1}^\ell$, $\left\{ \mu_i \right\}_{i=1}^\ell \subseteq P^+$ such that $\lambda_i - \mu_i \in Q$ for all $i$, and 
\[ \bpi_1 = \prod_{i=1}^\ell \bpi_{\lambda_i, a_i}, \quad \bpi_2 = \prod_{i=1}^\ell \bpi_{\mu_i, a_i} \ \in \P_{\text{Asym}}.\]
Then the $L^\sigma(\lie g)$--modules $V\left( \bor(\bpi_1) \right),   V\left( \bor(\bpi_2) \right)$ are strongly linked. 
\end{prop}
\begin{pf}
Again it is sufficient to prove the lemma for $\ell = 2$.  Since $\lambda_1 - \mu_1 \in Q$, by Proposition \ref{shaded}(ii) there exists a sequence $\left\{ \nu_i \right\}_{i = 0}^q \subseteq P^+$ with $\nu_0 = \lambda_1$, $\nu_q = \mu_1$ such that $\text{Hom}_{\lie g}(\lie g \otimes V(\nu_i), V(\nu_{i+1})) \neq 0$ for all $0 \leq i \leq q-1$. Fix some $i$.  By Lemma \ref{shaded}(i) either $\nu_i \geq \nu_{i+1}$ or $\nu_{i+1} \geq \nu_{i}$.  In either case, $V(\bor(\bpi_{\nu_i, a_1}))$ and  $V(\bor(\bpi_{\nu_{i+1}, a_1}))$ are both simple constituents of an indecomposable module $M_i$, which is in turn a quotient of $W(\bor(\bpi_{\nu_i, a_1}))$ for $\nu_i \geq \nu_{i+1}$ or $W(\bor(\bpi_{\nu_{i+1}, a_1}))$ for $\nu_i \leq \nu_{i+1}$.  We may assume without loss of generality that $\nu_i \geq \nu_{i+1}$.  Therefore $V(\bor(\bpi_{\nu_{i}, a_1})) \otimes V(\bor(\bpi_{\lambda_2,a_2}))$ and $V(\bor(\bpi_{\nu_{i+1}, a_1})) \otimes V(\bor(\bpi_{\lambda_2,a_2}))$ are both simple constituents of $M_i \otimes V(\bor(\bpi_{\lambda_2,a_2}))$, which is in turn a quotient of $W(\bor(\bpi_{\nu_i, a_1})) \otimes W(\bor(\bpi_{\lambda_2, a_2})) \cong W(\bor(\bpi_{\nu_i, a_1} \bpi_{\lambda_2, a_2}))$, hence indecomposable.  Therefore the modules 
\begin{align*} &V(\bor(\bpi_{\lambda_1, a_1})) \otimes V(\bor(\bpi_{\lambda_2, a_2})) \cong V(\bor(\bpi_{\lambda_1, a_1}\bpi_{\lambda_2, a_2})),\\
& \ \ \ \ \  \ \ V(\bor(\bpi_{\mu_1, a_1})) \otimes V(\bor(\bpi_{\lambda_2, a_2})) \cong V(\bor(\bpi_{\mu_1, a_1}\bpi_{\lambda_2, a_2})) \end{align*}
are strongly linked.  Similarly we can show that the modules 
\[  V(\bor(\bpi_{\mu_1, a_1}\bpi_{\lambda_2, a_2})), \ \  V(\bor(\bpi_{\mu_1, a_1}\bpi_{\mu_2, a_2})) \]
are strongly linked, completing the proof.
\end{pf}

\begin{cor}
Let $\bpi \in \P_{\text{Asym}}$, $\chi_{\bpi} = 0$.  Then $V(\bor(\bpi))$ is strongly linked to $\mathbb{C}$. 
\begin{pf}
The result follows from the above proposition \ref{mainprop} and two observations.  First, if $\chi_{\bpi} = 0$ then $\bpi$ is of the form $\prod \bpi_{\lambda_i, a_i}$, $\lambda_i \in Q \cap P^+$; and second, for $\lambda \in Q \cap P^+$, $a \in \mathbb{C}^*$, $V(\bor(\bpi_{\lambda, a}))$ is strongly linked to $\mathbb{C}$. 
\end{pf}
\end{cor}

\begin{cor}\label{irredlinkage}
Let $V(\bpi^\sigma_1), V(\bpi^\sigma_2) \in \cal F^\sigma_{\overline{\chi}}$.  then $V(\bpi^\sigma_1), V(\bpi^\sigma_2)$ are strongly linked.
\end{cor}
\begin{pf}
Let $\bpi_i \in \bor^{-1}(\bpi^\sigma_i)_{\text{Asym}}$, $i = 1,2$, so that we have 
$$
\xymatrix{
 & &  \bpi_1 \ar[dl]^{\bor} \ar@{<->}[r]^{\chi}  & \tilde \bpi_1 \ar[dr]_{\bor} & & \tilde \bpi_2 \ar[dl]^{\bor} \ar@{<->}[r]^{\chi} & \bpi_2 \ar[dr]_{\bor}  & \\
  & \bpi^\sigma_1  &  &  & \tilde{\bpi}^\sigma & & & \bpi^\sigma_2& 
}
$$
By Lemma \ref{sqarefreerep}, we can assume without loss of generality that $\tilde \bpi_1, \tilde \bpi_2 \in \P_{\text{Asym}}$.  It suffices now to show that $V(\bpi^\sigma_1)$ is strongly linked to $V(\tilde \bpi^\sigma)$. 

Let 
\[ \bpi^\sigma_1  = \prod_{i = 1}^\ell\prod_{\ep = 0}^{m-1}  \bpi^\sigma_{\lambda_{i, \ep}, \zeta^\ep a_i}, \ \ \ a_i^m \neq a_j^m \]
be a factorization of $\bpi^\sigma_1$.  Then we have $\bpi_1$ =$\prod_{i=1}^\ell \bpi_{\lambda_i, a_i}$ for some $\lambda_i \in P^+$, and $\tilde\bpi_1 =  \prod_{i=1}^\ell \bpi_{\lambda_i + \lambda_i', a_i} \tilde \bpi$, where $\lambda_i' \in Q$ such that $\lambda_i + \lambda_i' \in P^+$, $\chi_{\tilde\bpi} =0$ and the coordinates $\{b_i\}$ of $\tilde \bpi$ all satisfy $b^m \neq a^m$.  Furthermore 
\[ V(\tilde\bpi^\sigma) = V(\bor(\tilde \bpi_1)) = V(\bor(\prod_{i=1}^\ell \bpi_{\lambda_i + \lambda_i', a_i} \tilde \bpi)) \cong \bigotimes_{i = 1}^\ell  V(\bor(\bpi_{\lambda_i + \lambda_i', a_i})) \otimes V( \bor(\tilde \bpi)). \]
since $\chi_{\tilde\bpi} = 0$, we can conclude from Proposition \ref{mainprop} and the corollary following that the modules
\begin{align*}  
 \bigotimes_{i = 1}^\ell  V(\bor(\bpi_{\lambda_i + \lambda_i', a_i}))  \cong  & \bigotimes_{i = 1}^\ell  V(\bor(\bpi_{\lambda_i + \lambda_i', a_i})) \otimes \mathbb{C}, \\
& \bigotimes_{i = 1}^\ell  V(\bor(\bpi_{\lambda_i + \lambda_i', a_i}))\otimes V( \bor(\tilde \bpi)) \cong V(\tilde\bpi^\sigma)
\end{align*}
are strongly linked, as are the modules $V(\bpi^\sigma_1) \cong V(\bor(\bpi_1))$  and $\bigotimes_{i = 1}^\ell  V(\bor(\bpi_{\lambda_i + \lambda_i', a_i}))$.  This concludes the proof.
\end{pf}

\subsection{Proof of Proposition \ref{secondprop}} We first prove one important result concerning $\ext(U,V)$ for modules $U, V \in \cal F^\sigma$: that `distinct spectral characters have no non--trivial extensions' (Lemma \ref{extresults} (ii)).  We begin with a lemma required for the proof this fact.  In the following, $w_0$ is the longest element of the Weyl group of $\lie g$, and for a standard decomposition $\bpi = \prod_{i=1}^\ell \bpi_{\lambda_i, a_i}$, we define $\bpi^* = \prod_{i=1}^\ell \bpi_{-w_0 \lambda_i, a_i}$.  Then $V(\bpi)^* \cong V(\bpi^*)$ (\cite{CM}, prop. 3.2).  Also it is easy to see that $\lambda_\bpi = \sum_{i=1}^\ell \lambda_i$.   For any irreducible $L^\sigma(\lie g)$--module $V(\bpi^\sigma)$, let $(\bpi^\sigma)^*$ be the element of $\P^\sigma$ such that $V((\bpi^\sigma)^*) \cong V(\bpi^\sigma)^*$.

\begin{lem}\label{dual}\mbox{}
\begin{enumerate}
\item[(i)] $V(\bor(\bpi))^* \cong V(\bor(\bpi^*))$
\item[(ii)] $\displaystyle{\lambda_{(\bpi^\sigma)^*} = \lambda_{\bpi^\sigma}}$
\end{enumerate}
\end{lem}
\begin{pf}  For any $\bpi \in \P$, we have $V(\bpi)^* \cong V(\bpi^*)$ (\cite{CM}, prop. 3.2).  Therefore 
\[ V(\bor(\bpi))^* \cong \left( V(\bpi)\left. \right|_{L^\sigma(\lie g)}\right)^* \cong \left( V(\bpi)^* \right)\left. \right|_{L^\sigma(\lie g)} \cong V(\bpi^*)\left. \right|_{L^\sigma(\lie g)} \cong V(\bor(\bpi^*)).\]  
For the proof of (ii), let $\bpi \in \bor^{-1}(\bpi^\sigma)$.  Then $\lambda_{\bpi^\sigma} = \sum_{i=0}^\epsilon \lambda_{\bpi}(\epsilon)$.  For $\lie g$ of type $A$, $D$ or $E_6$, we have either $-w_0 = Id$ or $-w_0 = \sigma$ (see \cite{Bour}, for ex).  In either case, for any $\lambda \in P^+$ we have
\[ \sum_{\epsilon = 0}^{m-1}-w_0\lambda(\epsilon) = \sum_{\epsilon = 0}^{m-1}\lambda(\epsilon).\] 
 Also for any $\lambda, \mu \in P^+$, we have $(\lambda + \mu)(\ep) = \lambda(\ep) + \mu(\ep)$, $0 \leq \ep \leq m-1$.  Therefore
\[ \lambda_{(\bpi^\sigma)^*} = \lambda_{\bor(\bpi^*)} =  \sum_{\epsilon=0}^{m-1} \lambda_{\bpi^*}(\epsilon)  =  \sum_{\epsilon=0}^{m-1}  \lambda_{\bpi}(\epsilon) = \lambda_{\bpi^\sigma},
\]
where for the first equality we have used (i) of the lemma.
\end{pf}

\noindent The proof of the following is an adaptation of that given in \cite{CM}. 
\begin{lem}\label{extresults}
\mbox{}
\begin{enumerate}
\item [(i)] Let $U \in \cal F^\sigma_{\overline{\chi}}$, and $\bpi^\sigma \in \P^\sigma$ such that $\overline{\chi} \neq \chi_{\bpi^\sigma}$.  Then $\operatorname{Ext}^1_{L^\sigma(\lie g)}(U, V(\bpi^\sigma))=0$. 
\item [(ii)] Assume that $V_j \in \cal F^\sigma_{\overline{\chi_j}}$, $j = 1,2$ and that $\overline{\chi_1} \neq \overline{\chi_2}$.  Then $\operatorname{Ext}^1_{L^\sigma(\lie g)}(V_1, V_2) =0$.
\end{enumerate}
\end{lem}
\begin{pf}
Since $\operatorname{Ext}^1$ preserves direct sums, to prove the lemma it suffices to consider the case when $U$ is indecomposable.  Consider an extension 
\[ 0 \longrightarrow V(\bpi^\sigma_1) \longrightarrow V \longrightarrow U \longrightarrow 0.\]
We prove by induction on the length of $U$ that this extension must be trivial.  So first suppose that $U = V(\bpi^\sigma_2)$ for some $\bpi^\sigma_2 \in \P^\sigma$ and that $\chi_{\bpi^\sigma_2} \neq \chi_{\bpi^\sigma_1}$, so we have 
\[ 0 \longrightarrow V(\bpi^\sigma_1) \stackrel{\iota}{\longrightarrow} V \stackrel{\pi}{\longrightarrow} V(\bpi^\sigma_2) \longrightarrow 0.\]
For the remainder of the proof, let $\lambda_i = \lambda_{\bpi^\sigma_i} \in P_0^+$.  We must have either 
\begin{enumerate}
\item $\lambda_2  < \lambda_1$, or 
\item $\lambda_1  - \lambda_2 \notin (Q_0^+ - \left\{ 0 \right\})$.
\end{enumerate}
If we are in case (1) then dualizing the above exact sequence takes us to 

\[ 0 \longrightarrow V(\bpi^\sigma_2)^* \longrightarrow V^* \longrightarrow V(\bpi^\sigma_1)^* \longrightarrow 0\]
which, by Lemma \ref{dual}, takes us to case (2), so we can assume without loss of generality that we are in case (2). The exact sequence always splits as a sequence of $\lie g_0$-modules, so we have 
\[ V \cong_{\lie g_0} V(\bpi^\sigma_1)_{\lie g_0} \oplus V(\bpi^\sigma_2)_{\lie g_0} .\]
Therefore 
\[ V_{\lambda_2  } \cong V(\bpi^\sigma_1)_{\lambda_2  }  \oplus V(\bpi^\sigma_2)_{\lambda_2  }  .\]
Since we are in case (2), we know that $\lambda_2   \notin wt(V(\bpi^\sigma_1))$, and therefore
\[ V_{\lambda_2 } \cong V(\bpi^\sigma_2)_{\lambda_2},\]
and hence $L^\sigma(\lie n^+)V_{\lambda_2} = 0$.  On the other hand, since $V_{\lambda_2  }$ maps onto $V(\bpi^\sigma_2)_{\lambda_2 }$, there must be some nonzero vector $v \in V_{\lambda_2}$ with $L^\sigma(\lie h)$-eigenvalue $\bpi^\sigma_2$.  Therefore the submodule $\bu(L^\sigma(\lie g)).v$ of $V$ must be a quotient of $W(\bpi^\sigma_2)$, hence $\bu(L^\sigma(\lie g)).v \in \cal F^\sigma_{\chi_{\bpi^\sigma_2}}$.  If $\bu(L^\sigma(\lie g)).v=V$, then $V$ has spectral character $\chi_{\bpi^\sigma_2}$, but $V(\bpi^\sigma_1)$ is a submodule of $V$ and $\chi_{\bpi^\sigma_2} \neq \chi_{\bpi^\sigma_1}$.  Therefore $\bu(L^\sigma(\lie g)).v$ must be a proper nontrivial submodule of $V$.  But then $l(V)=2$ implies that either $\bu(L^\sigma(\lie g)).v \cong V(\bpi^\sigma_1)$ or $\bu(L^\sigma(\lie g)).v \cong V(\bpi^\sigma_2)$, and since $\chi_{\bpi^\sigma_2} \neq \chi_{\bpi^\sigma_1}$ we must have $\bu(L^\sigma(\lie g)).v \cong V(\bpi^\sigma_2)$.  Also we have  $\iota(V(\bpi^\sigma_1)) \cap \bu(L^\sigma(\lie g)).v = 0$, hence $V \cong V(\bpi^\sigma_1) \oplus V(\bpi^\sigma_2)$, and the induction begins.

Now assume that $U$ is indecomposable of length $\geq 1$ and $U \in \cal F^\sigma_{\overline{\chi}}$.  Let $U_1$ be a proper non-trivial submodule of $U$ and consider the short exact sequence 
\[ 0 \longrightarrow U_1 \longrightarrow U \longrightarrow U_2 \longrightarrow 0,\]
where $U_2 = U/U_1$. Since $U \in \cal F^\sigma_{\overline{\chi}}$, $U_i \in \cal F^\sigma_{\overline{\chi}}$ as well.  Then the inductive hypothesis gives us $\operatorname{Ext}^1_{L^\sigma(\lie g)}(U_i, Vr(\bpi^\sigma_1)) = 0$, and the result follows by using the exact sequence 
\[0 \longrightarrow  \operatorname{Ext}^1_{L^\sigma(\lie g)}(U_2, V(\bpi^\sigma_1)) \longrightarrow \operatorname{Ext}^1_{L^\sigma(\lie g)}(U, V(\bpi^\sigma_1)) \longrightarrow \operatorname{Ext}^1_{L^\sigma(\lie g)}(U_1, V(\bpi^\sigma_1)) \longrightarrow 0 .\]
 Part (ii) is now immediate by using a similar induction on the length of $V_2$. 
\end{pf}

\noindent We now conclude with the proof of Proposition \ref{secondprop}: let $V$ be an indecomposable $L^\sigma(\lie g)$-module.  We will show that there exists $\chi \in \Xi$ such that $V \in \cal F^\sigma_{\overline{\chi}}$ by an induction on the length of $V$.  If $V$ is irreducible, the result is immediate.  Now assume $V$ is reducible, and let $V(\bpi^\sigma)$ be an irreducible submodule of $V$; let $U = V/V(\bpi^\sigma)$.  So we have an extension 

\[ 0 \longrightarrow V(\bpi^\sigma) \longrightarrow V \longrightarrow U \longrightarrow 0.\]
Now decompose $U$:
\[ U = \bigoplus_{j=1}^rU_j; \ \ U_j \text{ indecomposable}.\]
Clearly $l(U_j) < l(V)$; therefore the inductive hypothesis ensures that $U_j \in \cal F^\sigma_{\overline{\chi_j}}$ for some $\chi_j \in \Xi$; $1 \leq j \leq r$.  Now we would like to argue that $\overline{\chi_j} = \chi_{\bpi^\sigma}$ for all $j$, for if so, then $U_j \in \cal F^\sigma_{\chi_{\bpi^\sigma}}$ for all $j$ and hence $U \in \cal F^\sigma_{\chi_{\bpi^\sigma}}$. 

Suppose instead that there exists some $j_0$ such that $\overline{\chi_{j_0}} \neq \chi_{\bpi^\sigma}$.  Then Lemma \ref{extresults} gives us
$$\operatorname{Ext}^1_{L^\sigma(\lie g)}(U,V(\bpi^\sigma)) \cong \bigoplus_{j=1}^r 
\operatorname{Ext}^1_{L^\sigma(\lie g)}(U_j,V(\bpi^\sigma))\cong
\bigoplus_{j\neq j_0} \operatorname{Ext}^1_{L^\sigma(\lie g)}(U_j,V(\bpi^\sigma));$$
I.e., the exact sequence $0\to V(\bpi^\sigma)\to V\to U \to 0$ is
equivalent to one of the form
\begin{equation*}
0\to V(\bpi^\sigma)\to U_{j_0}\oplus V'\to U_{j_0}\bigoplus_{j\neq j_0} U_j \to 0
\end{equation*}
where
\begin{equation*}
0\to V(\bpi^\sigma)\to V'\to \bigoplus_{j\neq j_0} U_j \to 0
\end{equation*}
is an element of $\oplus_{j\neq j_0} \operatorname{Ext}^1_{L^\sigma(\lie g)}(U_j,V(\bpi^\sigma))$. 
But this contradicts
the indecomposability of $V$. Hence $\overline{\chi_j}=\chi_{\bpi^\sigma}$ for all 
$1\le j\le r$ and $V\in\cal F^\sigma_{\chi_{\bpi^\sigma}}$.

\end{document}